\documentclass[11pt]{article}
\usepackage{amsmath,amsfonts,amssymb}
\usepackage{fancyhdr,graphicx,amsmath, amscd}
\usepackage{color}

\def\be{\begin{equation}}
\def\ee{\end{equation}}
\def\bea{\begin{eqnarray}}
\def\eea{\end{eqnarray}}
\def\bt{\begin{theorem}}
\def\et{\end{theorem}}
\def\bl{\begin{lemma}}
\def\el{\end{lemma}}
\def\br{\begin{remark}}
\def\er{\end{remark}}
\def\bc{\begin{corollary}}
\def\ec{\end{corollary}}
\def\bd{\begin{definition}}
\def\ed{\end{definition}}
\def\a{\alpha}
\def\b{\beta}
\def\g{\gamma}
\def\d{\delta}

\def\k{\kappa}
\def\l{\lambda}
\def\m{\mu}

\def\r{\rho}
\def\s{\sigma}

\def\D{\Delta}

\def\cE{\mathcal{E}}

\def\bbN{\mathbb{N}}

\def\bbR{\mathbb{R}}

\def\KC{\partial K \cap
\overline{[x,K]\backslash K}}

\def\bn{B_2^}

\def\b1{B_{1}^}

\def\pt{\partial}

\def\ba{\begin{array}}
\def\ea{\end{array}}
\def\ben{\begin{enumerate}}
\def\een{\end{enumerate}}
\newtheorem{theorem}{Theorem}[section]
\newtheorem{lemma}{Lemma}[section]
\newtheorem{remark}{Remark}[section]
\newtheorem{proposition}{Proposition}[section]
\newtheorem{corollary}{Corollary}[section]
\newtheorem{example}{Example}[section]
\newtheorem{definition}{Definition}[section]
\begin{document}
\thispagestyle{fancy}
\title{Inequalities for mixed $p$-affine surface area
}

\author{Elisabeth Werner\thanks{Partially supported by an NSF grant, a FRG-NSF grant and  a BSF grant}
\and Deping Ye }
\date{}

\maketitle
\begin{abstract}
We prove new Alexandrov-Fenchel type inequalities and new affine
isoperimetric inequalities for mixed $p$-affine surface areas. We
introduce a new class of bodies, the illumination surface bodies,
and establish some of their properties. We show, for instance, that they are not
necessarily convex. We give geometric interpretations of $L_p$ affine surface areas, mixed
$p$-affine surface areas and  other functionals   via these
bodies. The surprising new element is that not necessarily convex
bodies provide the tool for these interpretations.
\end{abstract}

\section{Introduction}

This article deals with affine isoperimetric inequalities and
Alexandrov-Fenchel type inequalities for mixed $p$-affine surface
area. Mixed $p$-affine surface area was introduced by Lutwak  for
$p \geq 1$ in  \cite{Lu1}. It has the dual mixed volume
\cite{Lut1975} and the $L_p$ affine surface area  \cite{Lu1} as
special cases. $L_p$ affine surface area is at the core of the
rapidly developing $L_p$ Brunn-Minkowski theory. Contributions
here include the study of solutions of nontrivial ordinary and,
respectively, partial differential equations (see e.g. Chen
\cite{CW1}, Chou and Wang \cite{CW2}, Stancu \cite{SA1, SA2}), the
study of the $L_p$ Christoffel-Minkowski problem by Hu, Ma and
Shen \cite{HMS}, extensions of $L_p$ affine surface area to all
$p$  (see e.g., \cite{MW2, SW4, SW5,
WY}), a new proof by Fleury, Gu\'edon and Paouris \cite{FGP}  of a
result by Klartag \cite{Kl} on concentration of volume,  results
on approximation of convex bodies by polytopes (e.g., \cite{ Gr2,
LuSchW, SW5}), results on valuations (e.g., Alesker  \cite{A1, A2}, and
Ludwig and Reitzner \cite{LudR, LR1})
and the affine Plateau problem
solved in ${\mathbb R}^3$ by Trudinger and Wang \cite{TW1}, and
Wang \cite{Wa}.

\par The classical affine isoperimetric inequality, which gives an
upper bound for the  affine surface area in terms
of volume, is fundamental in many problems (e.g. \cite{Ga, GaZ,
Lu-O, Sch}). In particular, it was used to show the uniqueness of
self-similar solutions of the affine curvature flow and to study
its asymptotic behavior by Andrews \cite{AN1, AN2}, Sapiro and
Tannenbaum \cite{ST1}. More general $L_p$ affine isoperimetric
inequalities were proved in \cite{Lu1} for $p>1$ and in \cite{WY}
for all $p$. These $L_p$ affine isoperimetric inequalities
generalize the celebrated Blaschke-Santal\'{o} inequality
and inverse Santal\'{o}
inequality due to Bourgain and Milman \cite{BM} (see also
Kuperberg \cite{GK2}). We also refer to related works by Lutwak,
Yang and Zhang \cite{LYZ1} and Campi and Gronchi \cite{CG}.

\par
For mixed $p$-affine surface area, Alexandrov-Fenchel type
inequalities (for $p=1, \pm \infty$) and affine isoperimetric
inequalities (for $1 \leq p \leq n$) were first established by
Lutwak in \cite{Lut1975, Lut1987, Lu1}. Here we derive new
Alexandrov-Fenchel type inequalities for mixed $p$-affine surface
area for all $p \in [-\infty, \infty]$ and new mixed  $p$-affine
isoperimetric inequalities for all $p\in [0, \infty]$.
Classification of the equality cases for all $p$  in the
Alexandrov-Fenchel type inequalities for mixed $p$-affine surface
area is related to the uniqueness of solutions of the $L_p$
Minkowski problem (e.g.,  \cite{CW1, CW2, Lut1993, Lu-O, LYZ2,
LYZ2004, SA1, SA2}), which is unsolved for many cases. This is
similar to the classical Alexandrov-Fenchel inequalities for mixed
volume, where the complete classification of the equality cases is
also an unsolved problem.

\par  We also give new geometric interpretations for
functionals on convex bodies.  In particular, for $L_p$ affine
surface area, mixed $p$-affine surface area, and $i$-th mixed
$p$-affine surface area (see below for the  definitions). To do
so, we construct a new class of bodies, the illumination surface
bodies, and study the asymptotic behavior of their volumes. We
show  that the illumination surface bodies are not necessarily
convex, thus introducing a novel idea in  the theory of geometric
characterizations of functionals on convex bodies, where to date
only convex bodies where used (e.g. \cite{MW2, SW4, SW5, WY}).

\vskip 3mm From now on, we will always assume that the centroid of
a convex body $K$ in $\mathbb R^n$ is at the origin. We write $K
\in C^2_+$,  if $K$ has $C^2$ boundary with everywhere strictly
positive Gaussian curvature. For real $p \geq 1$, the mixed
$p$-affine surface area, $as_{p}(K_1, \cdots,K_{n})$, of $n$
convex bodies  $K_i \in C^2_+$ was introduced in  \cite{Lu1} by
\begin{equation}\label{def:pmixed}
as_p(K_1, \cdots, K_{n})=\int
_{S^{n-1}}\bigg[h_{K_1}(u)^{1-p}f_{K_1}(u)\cdots
h_{K_n}^{1-p}f_{K_n}(u)\bigg]^{\frac{1}{n+p}}\,d\s (u).
\end{equation}
Here $S^{n-1}$ is the boundary of the Euclidean unit ball $B^n_2$
in  $\mathbb{R}^n$, $\s$ is the usual surface area measure on
$S^{n-1}$, $h_{K}(u)$ is the support function of the convex body
$K$ at $u\in S^{n-1}$, and $f_{K}(u)$ is the curvature function of
$K$ at
 $u$, i.e., the reciprocal of the Gauss
curvature $\kappa _K(x)$ at this point $x \in
\partial K$, the boundary of $K$,  that has $u$ as its outer normal.
\par
We propose here to extend the definition (\ref{def:pmixed}) for
mixed $p$-affine surface area to all $p\neq -n$. We also  propose
a definition for the $(-n)$- mixed affine surface area (see
Section 2).
\par
We show that  mixed $p$-affine surface
areas are affine invariants for all $p$.
Note  that for $p=\pm \infty$,
\begin{eqnarray}
as_{\pm \infty}(K_1, \cdots, K_{n})&=&\int _{S^{n-1}}
\frac{1}{h_{K_1}(u)} \cdots \frac{1}{h_{K_n}(u)} \,d\s (u)\nonumber \\
&=& n \tilde{V}(K_1^\circ, \cdots, K_n^\circ) \label{dual-mixed}
\end{eqnarray} where $K^\circ=\{y\in \bbR^n, \langle x, y\rangle\leq 1, \forall x\in
K\}$ is the polar body of $K$, and $\tilde{V}(K_1^\circ, \cdots,
K_n^\circ)$ is the dual mixed volume of $K_1^\circ, \cdots,
K_n^\circ$,  introduced by Lutwak in \cite{Lut1975}.

\vskip2mm \par
 When all $K_i$ coincide with $K$, then for all $p\neq -n$
\begin{equation}\label{all equal}
as_p(K, \cdots, K)=\int_{S^{n-1}}\frac{f_{K}(u)^{\frac{n}{n+p}}}
{h_K(u)^{\frac{n(p-1)}{n+p}}} d\sigma(u) =as_p(K).
\end{equation}
 $as_p(K)$ is the $L_p$ affine surface area of $K$, which is
defined for a general convex body $K$ as in \cite{Lu1} ($p >1$) and
in \cite{SW5} ($p <1$) by
\begin{equation} \label{def:paffine}
as_{p}(K)=\int_{\partial K}\frac{\kappa_K(x)^{\frac{p}{n+p}}}
{\langle x,N_{ K}(x)\rangle ^{\frac{n(p-1)}{n+p}}} d\mu_{ K}(x).
\end{equation}
$N_K(x)$ is the outer unit normal vector at $x$ to $\partial K$,
$\m _K$ denotes the usual surface area measure on $\pt K$, and
$\langle \cdot, \cdot\rangle$ is the standard inner product on
$\bbR^n$ which induces the Euclidian norm $\|\cdot\|$. If $K\in
C^2_+$, (\ref{def:paffine}) can be rewritten as (\ref{all equal}).
We show in Section 2 that the corresponding formula (\ref{all
equal}) for $p =-n$ also holds, where $as_{-n}(K)$ is the $L_{-n}$
 affine surface area of $K$ introduced in  \cite{MW2}.

\par Note further that the surface area of $K$ can be written as $(-1)$-th mixed
$1$-affine surface area of $K$ and the Euclidean ball $B^n_2$ (see
Section 2).
\par
Thus, mixed $p$-affine surface area is an
extension of dual mixed volume, surface area, and $L_p$ affine
surface area.

\vskip 3mm {\bf Further notations.} For sets $A$ and $B$,
$[A,B]=\mbox{conv}(A,B):=\{\l x+(1-\l)y: \l\in[0,1], x,y \in A\cup
B \}$ is the convex hull of $A\cup B$. A subset $K$ of $\bbR^n$ is
star convex if there exists $x_0\in K$ such that the line segment
$[x_0,x]$, from $x_0$ to any point $x$ in $K$, is contained in
$K$. A convex body $K$ is said to be strictly convex if $\pt K$
does not contain any line segment.
\par
For a convex  body $K$ in $\mathbb{R}^n$, $|K|$ stands for the $n$-dimensional
volume of $K$. More generally, for a set $M$, $|M|$
denotes the Hausdorff content of its appropriate dimension.

\par
For $u\in S^{n-1}$, $H(x,u)$ is the hyperplane  through $x$ with
outer normal vector $u$,  $H(x,u)=\{y\in \bbR^n, \langle y, u
\rangle=\langle x,u\rangle\}$. The two half-spaces generated by
$H(x,u)$ are $H^-(x,u)=\{y\in \bbR^n, \langle y, u \rangle\geq
\langle x,u\rangle\}$ and $H^+(x,u)=\{y\in \bbR^n, \langle y, u
\rangle\leq \langle x,u\rangle\}$. For $f: \partial K \rightarrow
\mathbb{R}_+ \cup \{0\}$, $\mu_f$ is the measure on $\partial K$
defined by $\m_f(A)= \int_A f d\mu_K$.

\vskip 3mm \par The paper is organized as follows. In Section 2,
we prove new  Alexandrov-Fenchel type  inequalities and new
isoperimetric inequalities for  mixed $p$-affine surface areas. We
show monotonicity behaviour of the quotients
$$\left(\frac{as_p(K_1, \cdots, K_n)}{as_\infty(K_1, \cdots, K_n)}\right)^{n+p}
\mbox{ and } \left(\frac{as_p(K_1, \cdots, K_n)}{as_0(K_1, \cdots,
K_n)}\right)^\frac{n+p}{p} .$$ We prove Blaschke-Santal\'{o} type
inequalities for  mixed $p$-affine surface areas. Similar results
for the $i$-th mixed $p$-affine surface areas are also proved in
 Section 2. In Section 3, we introduce the illumination surface body
 and describe some of its properties.
In Section 4, we derive the asymptotic behavior of the volume of
the illumination surface body, and geometric interpretations of
$L_p$ affine surface areas, mixed $p$-affine surface areas,  and
other functionals on convex bodies.

\section{ Mixed $p$-affine surface area and related inequalities}
\subsection{Inequalities for mixed $p$-affine surface area}

We begin by proving that mixed $p$-affine surface area is affine
invariant for all $p$. For $p \geq 1$, this was proved by Lutwak
\cite{Lu1}. We will first treat the case $p \neq -n$.  All the
results concerning the  case $p=-n$ are at the end of this
subsection.
\par
\noindent
It will be convenient to use the notation
\begin{equation}\label{fp}
 f_p(K,u)=h_K^{1-p}(u)f_K(u)
\end{equation}
for a convex body $K$ in $\mathbb{R}^n$ and  $u \in S^{n-1}$. We
will also write  $as_p^m(K_1, \cdots, K_n)$  for $\big[as_p(K_1,
\cdots, K_n)\big]^m$, and $|det(T)|$ for the absolute value
of the determinant of linear transform $T$.

\bl Let $T: \bbR^n \rightarrow\bbR^n $ be an invertible linear
transform. Then for all $p\neq -n$, $$as_p(TK_1, \cdots,
TK_{n})=|det(T)|^{\frac{n-p}{n+p}}\  as_p(K_1, \cdots,K_{n}).$$ In
particular, if $|det(T)|=1$, then $as_p(K_1, \cdots, K_n)$ is
affine invariant:
$$ as_p(TK_1, \cdots, TK_{n})= as_p(K_1,
\cdots,K_{n}).$$ \el

\vskip 2mm
\noindent
{\bf Proof.}
\par
\noindent
Since $K\in C^2_+$, for any $u\in S^{n-1}$, there
exists a unique $x\in \partial K$ such that $u=N_K(x)$ and
$f_K(u)=\frac{1}{\k _K(x)}$. By Lemma 12 of \cite{SW5}
\begin{eqnarray}\label{Curvature:K:TK}
f_K(u)=\frac{1}{\k_K(x)}= \frac{f_{TK}\left(v\right)}{\  det ^2(T)
\ \|T^{-1t}(u)\|^{n+1} },
\end{eqnarray} where $v=\frac{T^{-1t}(u)}{\|T^{-1t}(u)\|}\in S^{n-1}$ and
where for an operator $A$, $A^t$ denotes its usual adjoint.  On
the other hand,
$$h_K(u)=\langle x, u \rangle =\langle Tx, T^{-1t}(u)\rangle=\|T^{-1t}(u)\|\ \langle Tx, v \rangle=\|T^{-1t}(u)\|\ h_{TK}(v).$$

Thus, with  notation (\ref{fp}),  for all $p$,
\begin{equation}\label{Curvature:K:TK:2}f_p(K,u)=
\frac{f_{TK}\left(v\right)\ h_{TK}^{1-p}(v)\
\|T^{-1t}(u)\|^{1-p}}{det ^2(T) \ \|T^{-1t}(u)\|^{n+1}
}=\frac{f_p(TK,v)}{det^2(T)\ \|T^{-1t}(u)\|^{n+p}}.
\end{equation}
Lemma 10 and its proof in \cite{SW5} show that -up to a small
error-
$$f_{TK}(v)\,d\s(v)=|{det}(T)|\ \|T^{-1t}(u)\|f_K(u)\,d\s(u).$$ Together
with (\ref{Curvature:K:TK}), one gets that (again up to a small
error) $ \|T^{-1t}(u)\|^{-n}\,d\s(u)=|{det}(T)|\ \,d\s(v).$
Therefore, up to a small error,
\begin{eqnarray*} [f_p(K_1, u)\cdots
f_p(K_n,u)]^{\frac{1}{n+p}}\,d\s (u)=|{det}(T)|^{\frac{p-n}{n+p}}\
[f_p(TK_1, v)\cdots f_p(TK_n,v)]^{\frac{1}{n+p}}\,d\s (v).
\end{eqnarray*} The lemma then follows by  integrating  over $S^{n-1}$.

\vskip 3mm A general version of the classical Alexandrov-Fenchel
inequalities for  mixed volumes (see \cite{ Ale1937, Bus1958,
Sch}) can be written as
$$
\prod_{i=0}^{m-1}V(K_1,\cdots, K_{n-m}, \underbrace{K_{n-i},
\cdots, K_{n-i}}_m)\leq V^m(K_1, \cdots, K_n). $$ Here we prove
the analogous inequalities for mixed $p$-affine surface area. For
$p=\pm \infty$ and $p=1$, the  inequalities were proved by Lutwak
\cite{Lut1975, Lut1987}. For $p\geq 1$, inequality
(\ref{Holder:pmixed}) was proved by Lutwak in \cite{Lu1}, with
equality if and only if the associated $K_i$ are dilates of each
other.

\begin{proposition} \label{ineq1:pmixed}
Let all $K_i$ be convex bodies  in $C^2_+$ with centroid  at the
origin. If $p\neq -n$, then for $1 \leq m \leq n$
$$ as_p^m(K_1, \cdots, K_{n}) \leq \prod_{i=0}^
{m-1}as_p(K_1, \cdots, K_{n-m}, \underbrace{K_{n-i}, \cdots,
K_{n-i}}_m). $$ Equality holds if the $K_k$, for $k=n-m+1, \cdots,
n$ are dilates of each other. If $m=1$, equality holds trivially.
\par \noindent In particular, if $m=n$,
\begin{eqnarray}\label{Holder:pmixed} as_p^n(K_1,
\cdots,K_{n}) \leq  as_p(K_1) \cdots as_p(K_n).
\end{eqnarray}
\par \noindent
\end{proposition}
 {\bf Proof.} Put $g_0(u)=[f_p(K_1,u)\cdots
f_p(K_{n-m},u)]^{\frac{1}{n+p}}$ and
for $i=0, \cdots, m-1$, put
$g_{i+1}(u)=[
f_p(K_{n-i},u)]^{\frac{1}{n+p}}$.  By
H\"{o}lder's  inequality (see \cite{HLP})
\begin{eqnarray*}
as_p(K_1, \cdots, K_n)&=& \int _{S^{n-1}} g_0(u) g_1(u) \cdots
g_{m}(u)\,d\s(u)\\ &\leq & \prod _{i=0}^{m-1} \left(\int
_{S^{n-1}} g_0(u) g_{i+1}^m(u)\,d\s(u)\right)^{\frac{1}{m}}\\ &=&
\prod_{i=0}^ {m-1} as_p^{\frac{1}{m}} (K_1, \cdots,K_{n-m},
\underbrace{K_{n-i} \cdots K_{n-i}}_m).
\end{eqnarray*}
As  $K_i \in C^2_+$, $f_p(K_{i},u)>0$ for all $i$ and all $u\in
S^{n-1}$. Therefore, equality in H\"{o}lder's  inequality holds if
and only if $g_0(u) g_{i+1}^m(u)=\l ^m g_0(u) g_{j+1}^m(u)$ for
some $\l>0$ and all $0\leq i\neq j\leq m-1$. This is equivalent to
$h_{K_{n-i}}(u)^{1-p}f_{K_{n-i}}(u)=\l
h_{K_{n-j}}(u)^{1-p}f_{K_{n-j}}(u)$ for all $0\leq i\neq j\leq
m-1$. This condition holds true if the $K_k$,  for $k=n-m+1,
\cdots, n$ are dilates of each other. \vskip 2mm \noindent {\bf
Remark.} It is an unsolved problem for many $p$ whether
$f_p(K,u)=\l f_p(L,u)$  guarantees that $K$ and $L$ are dilates of
each other. This is equivalent to the uniqueness of the solution
of the $L_p$ Minkowski problem: {\it for  fixed $\alpha \in \bbR$,
under which conditions on a continuous function $\g:
S^{n-1}\rightarrow (0, \infty)$, there exists a (unique) convex
body $K$ such that $h_K(u)^\a f_K(u)= \g(u)$ for all $u\in
S^{n-1}$.}  In many cases, the uniqueness of the solution is an
open problem. We refer to e.g., \cite{CW2, Lut1993, Lu-O, LYZ2004,
SA1, SA2} for detailed information and more references on the
subject. For $p \geq 1, p \neq n$,  the solution to the $L_p$
Minkowski problem is known to be unique and for $p=n$, the
solution is unique modulo dilates \cite{Lut1993}. Therefore, we
have the characterization of equality in Proposition
\ref{ineq1:pmixed} for $p \geq 1$.

\vskip 3mm
\noindent
Next,  we prove  affine isoperimetric inequalities
for mixed $p$-affine surface areas.
\begin{proposition}\label{Iso:pmixed:2}
Let all $K_i$ be convex bodies in $C^2_+$ with centroid at the
origin. \vskip 2mm
 \noindent
(i) For $ p\geq 0$,
\begin{eqnarray*}
\frac{as_p^n(K_1, \cdots, K_n)}{as_p^n(B^n_2, \cdots, B^n_2)}
\leq \left( \frac{|K_1|}{|B^n_2|}
\cdots \frac{|K_n|}{|B^n_2|}\right)^\frac{n-p}{n+p},
\end{eqnarray*}
with equality if the $K_i$ are
ellipsoids that are dilates of one another.
\vskip 2mm
 \noindent
 (ii) For $0\leq p\leq n$,
$$\frac{as_p(K_1, \cdots, K_n)}{as_p(B^n_2, \cdots, B^n_2)}\leq
\left(\frac{V(K_1, \cdots, K_n)}{V(B^n_2, \cdots,
B^n_2)}\right)^{\frac{n-p}{n+p}}$$ with equality if
the $K_i$ are ellipsoids that are dilates of one another.

\vskip 3mm \par \noindent In particular, for $p=n$
$$
as_n(K_1, \cdots, K_n) \leq as_n(B^n_2, \cdots, B^n_2),
$$
with equality if and only if the $K_i$ are ellipsoids that are
dilates of one another. \vskip 2mm \noindent (iii) For $p\geq n$,
$$
\frac{as_p(K_1,
\cdots, K_n)}{as_p(B^n_2, \cdots, B^n_2)}\leq
\left(\frac{\tilde{V}(K_1, \cdots, K_n)}{\tilde{V}(B^n_2, \cdots,
B^n_2)}\right)^{\frac{n-p}{n+p}},
$$ with equality if and only if
the $K_i$ are ellipsoids that are dilates of one another.
\vskip2mm \par \noindent In particular, for $p=\pm \infty$
$$ \tilde{V}(K_1, \cdots, K_n) \tilde{V}(K_1^\circ, \cdots, K_n^\circ) \leq |B_2^n|^2,
$$
with equality if and only if $K_i$ are ellipsoids that are dilates of one another.
\end{proposition}

\par \noindent {\bf Remark.} For $1\leq p\leq n$,
inequality (ii) (with equality if and only if the $K_i$ are
ellipsoids that are dilates of one another) was proved by Lutwak
in \cite{Lu1}. If $K_i=K$ for all $i$, one recovers the $L_p$
affine isoperimetric inequality proved in \cite{WY}.

\vskip 2mm \par \noindent {\bf Remark.} We  cannot expect to
get strictly positive lower bounds in Proposition
\ref{Iso:pmixed:2}. As  in \cite{WY}, we consider the
convex body $K(R, \varepsilon)\subset \bbR^2$, obtained as the
intersection of four Euclidean balls with radius $R$ centered at
$(\pm(R-1),0)$, $(0, \pm(R-1))$, $R$ arbitrarily large. We then
``round" the corners by putting there arcs of Euclidean balls of
radius
 $\varepsilon$, $\varepsilon$ arbitrarily small.
To obtain a body in $C^2_+$, we ``bridge" between the $R$-arcs and
$\varepsilon$-arcs  by $C^2_+$-arcs on a set of  arbitrarily small
measure. Then $as_p(K(R, \varepsilon)) \leq
\frac{16}{R^\frac{p}{2+p}} + 4  \pi\ \varepsilon^\frac{2}{2+p}$,
which goes to $0$ as $R\rightarrow \infty$ and $\varepsilon
\rightarrow 0$. Choose now $R_i$  and $\varepsilon_i$, $ 1 \leq i \leq n$,  such that
$R_i\rightarrow
\infty$ and $\varepsilon _i \rightarrow 0$,   and let
 $K_i= K(R_i,
\varepsilon _i) $ for $i=1, 2 \cdots, n$. By inequality
(\ref{Holder:pmixed}),  $as_p^n(K_1, \cdots, K_n)\leq  \prod
_{i=1}^n as_p(K(R_i, \varepsilon _i))$
and thus
$as_p(K_1, \cdots, K_n) \rightarrow 0$ for $p>0$. A similar
construction can be done in higher dimensions.

\vskip 3mm
\par
\noindent {\bf Proof of Proposition \ref{Iso:pmixed:2}.}
\newline
(i)
Clearly $as_p(B^n_2, \cdots,
B^n_2)=as_p(B^n_2)=n|B^n_2|$ for all $p\neq -n$. By inequality
(\ref{Holder:pmixed}), one gets for all $p \geq 0$
\begin{eqnarray}\label{Iso:pmixed:1}\frac{as_p^n(K_1, \cdots, K_n)}{as_p^n(B^n_2, \cdots, B^n_2)}
\leq \frac{as_p(K_1)}{as_p(B^n_2)} \cdots
\frac{as_p(K_n)}{as_p(B^n_2)}\leq \left( \frac{|K_1|}{|B^n_2|}
\cdots \frac{|K_n|}{|B^n_2|}\right)^\frac{n-p}{n+p}.
\end{eqnarray} The second inequality follows, for $p \geq 0$, from the $L_p$
affine isoperimetric inequality in \cite{WY}. Equality holds true
in the $L_p$ isoperimetric inequality \cite{WY} if and only if the
$K_i$ are all ellipsoids, and equality holds true in inequality
(\ref{Holder:pmixed}) if the $K_i$ are dilates of one another.
Thus, equality holds in (\ref{Iso:pmixed:1}) if the $K_i$ are
ellipsoids that are dilates of one another.

\vskip 2mm
\noindent
(ii) A direct consequence of the classical
Alexandrov-Fenchel inequality for mixed volume (see e.g.
\cite{Bur1988, Lei1980}) is that
$$|K_1|\cdots |K_n|\leq V^n(K_1, \cdots, K_n).$$
If $0\leq p \leq n$, then $\frac{n-p}{n+p}\geq 0$. Thus
$$
\big(|K_1|\cdots |K_n|\big)^\frac{n-p}{n+p} \leq \big[V^n(K_1,
\cdots, K_n)\big]^\frac{n-p}{n+p}.
$$
As $V(B^n_2, \cdots, B^n_2)=|B^n_2|$,  one
gets together with (\ref{Iso:pmixed:1})
$$
\frac{as_p(K_1, \cdots, K_n)}{as_p(B^n_2, \cdots, B^n_2)}\leq
\left(\frac{V(K_1, \cdots, K_n)}{V(B^n_2, \cdots,
B^n_2)}\right)^{\frac{n-p}{n+p}},
$$ with equality if the $K_i$ are
ellipsoids that are dilates of one another.
\vskip 2mm
\noindent
(iii) The analogous inequality for dual mixed volume
\cite{Lut1975} is
$$ |K_1|\cdots |K_n| \geq \tilde{V}^n(K_1,
\cdots, K_n),
$$
with equality if and only the $K_i$ are dilates of
one another.  $p > n$ implies $\frac{n-p}{n+p}< 0$. Thus
$$
\big(|K_1|\cdots |K_n|\big)^\frac{n-p}{n+p} \leq
\big[\tilde{V}^n(K_1, \cdots, K_n)\big]^\frac{n-p}{n+p}.$$
Together with (\ref{Iso:pmixed:1}) and $\tilde{V}(B^n_2, \cdots,
B^n_2)=|B^n_2|$, one gets $$\frac{as_p(K_1, \cdots,
K_n)}{as_p(B^n_2, \cdots, B^n_2)}\leq \left(\frac{\tilde{V}(K_1,
\cdots, K_n)}{\tilde{V}(B^n_2, \cdots,
B^n_2)}\right)^{\frac{n-p}{n+p}}.$$ As for $p \geq 1$ equality in
(\ref{Holder:pmixed}) holds if and only if the $K_i$ are dilates
of one another, equality holds true here if and only if the $K_i$
are ellipsoids that are dilates of one another.

\begin{proposition}
Let $\cE$ be a centered ellipsoid. If either all $K_i\in C^2_+$
are subsets of $\cE$, for $0\leq p<n$, or $\cE$ is subset of all
$K_i$ for $p>n$, then $$as_p(K_1, \cdots, K_n)\leq as_p(\cE).$$
For $p=n$, the inequality  holds for all $K_i$ in $C^2_+$ by Proposition \ref{Iso:pmixed:2} (ii).
\end{proposition}
\par
\noindent
{\bf Remark.} This proposition was proved by Lutwak  \cite{Lu1} if $K_i=K$ for
all $i$.
\vskip 2mm
\noindent
{\bf Proof.} It is
enough to prove the proposition for $\cE=B^n_2$. For $0\leq p<n$,
one has $\frac{n-p}{n+p}>0$ and hence
$\left(\frac{|K_i|}{|B^n_2|}\right)^{\frac{n-p}{n+p}}\leq 1$ as
$K_i\subset B^n_2$. Similarly, $p>n$ implies $\frac{n-p}{n+p}<0$
and therefore
$\left(\frac{|K_i|}{|B^n_2|}\right)^{\frac{n-p}{n+p}}\leq 1$ as
$B^n_2\subset K_i$ for all $i$. In both cases, the proposition
follows by inequality (\ref{Iso:pmixed:1}).

\vskip 3mm
\noindent
The next proposition gives a Blaschke-Santal\'{o} type inequality for  $p$-mixed affine surface area.
When $K_i=K$ for all $i$,
the proposition was proved in \cite{WY}.

\vskip 3mm
\noindent
\begin{proposition}\label{Santalo:0}
\par\noindent Let all $K_i$ be convex bodies  in $C^2_+$ with centroid
at the origin. For all $p\geq 0$,
\begin{eqnarray}\label{Santalo:1}
as_p^n(K_1, \cdots,
K_n)as_p^n(K_1^\circ, \cdots, K_n^\circ) \leq n^{2n}
|K_1||K_1^\circ|\cdots |K_n||K_n^\circ|.
\end{eqnarray}
Furthermore, $as_p(K_1, \cdots,
K_n)as_p(K_1^\circ, \cdots, K_n^\circ) \leq as_p^2(B^n_2)$ with
equality if the $K_i$ are ellipsoids that are dilates of one another.
\end{proposition}

\noindent  {\bf Proof.} It follows from  (\ref{Holder:pmixed})
that  for all $p\neq -n$,
$$as_p^n(K_1, \cdots, K_n)as_p^n(K_1^\circ, \cdots, K_n^\circ)\leq as_p(K_1)as_p(K_1^\circ)
\cdots as_p(K_n)as_p(K_n^\circ),
$$
with equality if the $K_i$ are
dilates of one another. By Corollary 4.1 in \cite{WY}, for $p\geq
0$,
$$as_p^n(K_1, \cdots, K_n)as_p^n(K_1^\circ, \cdots, K_n^\circ)\leq n^{2n} |K_1||K_1^\circ|\cdots |K_n||K_n^\circ|.$$
Blaschke-Santal\'{o} inequality states that
$|K||K^\circ|\leq |B^n_2|^2$ with equality if and only if  $K$ is a
$0$-centered ellipsoid. We apply it to inequality
(\ref{Santalo:1}), and obtain that for $p\geq 0$,
$$as_p(K_1, \cdots, K_n)as_p(K_1^\circ, \cdots, K_n^\circ) \leq
as_p^2(B^n_2, \cdots, B^n_2).$$ Equality holds if the $K_i$ are
ellipsoids that are dilates of one another.

\vskip 3mm
\noindent
\bt\label{ineq2:pmixed} Let $s\neq -n,  r \neq -n, p \neq -n$ be
real numbers. Let all $K_i$ be convex
bodies in $C^2_+ $ with centroid at the origin.
\par \noindent
(i)
If $\frac{(n+p)(r-s)}{(n+r)(p-s)}>1$, then
\begin{equation*}\label{i-2} as_p(K_1, \cdots,K_{n}) \leq \big(as_r(K_1, \cdots,K_{n})\big)^{\frac{(p-s)(n+r)}{(r-s)(n+p)}}
\big(as_s(K_1, \cdots,K_{n})\big)^{\frac{(r-p)(n+s)}{(r-s)(n+p)}}.
\end{equation*}
\par \noindent (ii) If $\frac{n+p}{n+r} >1$,  then
\begin{eqnarray}\label{i-3} as_p(K_1, \cdots,K_{n}) \leq \left(as_r(K_1, \cdots,
K_{n})\right)^{\frac{n+r}{n+p}} {\left(n\tilde{V}(K_1^\circ,
\cdots, K_{n}^\circ)\right)^{\frac{p-r}{n+p}}}. \end{eqnarray} \et
\vskip 2mm \noindent {\bf Remark.} When all $K_i$ coincide with
$K$, (i) of Theorem \ref{ineq2:pmixed} was proved in \cite{WY}.

\vskip 2mm \noindent {\bf Proof.} \vskip 2mm \noindent (i) By
H\"older's inequality -which enforces the condition
$\frac{(n+p)(r-s)}{(n+r)(p-s)}>1$,
\begin{eqnarray*}
as_p(K_1, \cdots,K_{n}) &=&\int _{S^{n-1}} \big[f_p(K_1,u) \cdots
f_p(K_n,u)\big]^\frac{1}{n+p}\,d\s (u)
\\ &=& \int _{S^{n-1}}
\left(\big[f_r(K_1,u) \cdots
f_r(K_n,u)\big]^\frac{1}{n+r}\right)^\frac{(n+r)(p-s)}{(n+p)(r-s)}\\
& & \ \ \left(\big[f_s(K_1,u) \cdots
f_s(K_n,u)\big]^\frac{1}{n+s}\right)^\frac{(n+s)(r-p)}{(n+p)(r-s)}
d\s (u) \\ &\leq & \big(as_r(K_1, \cdots, K_{n})
\big)^{\frac{(p-s)(n+r)}{(r-s)(n+p)}} \big(as_s(K_1, \cdots,
K_{n}) \big)^{\frac{(r-p)(n+s)}{(r-s)(n+p)}}.
\end{eqnarray*}

\par \noindent (ii) Similarly, again using H\"older's inequality -which now
enforces the condition $\frac{n+p}{n+r} >1$
\begin{eqnarray*}
as_p(K_1, \cdots, K_{n})&=&\int _{S^{n-1}} \big[f_p(K_1,u) \cdots
f_p(K_n,u)\big]^\frac{1}{n+p}\,d\s (u) \\
&=& \int _{S^{n-1}} \left(\big[f_r(K_1,u) \cdots
f_r(K_n,u)\big]^\frac{1}{n+r}\right)^\frac{n+r}{n+p} \\
& & \ \ \ \ \ \ \ \ \left[\frac{1}{h_{K_1}(u) \cdots h_{K_n}(u)} \right]^\frac{p-r}{n+p} \,d\s(u) \\
&\leq& \left(as_r(K_1, \cdots, K_{n})\right)^{\frac{n+r}{n+p}}\
{\left(as_{\infty}(K_1, \cdots, K_{n})\right)^{\frac{p-r}{n+p}}}.
\end{eqnarray*} Together with (\ref{dual-mixed}), this completes the proof.

\vskip 3mm \noindent{\bf Remark.} The condition
$\frac{(n+p)(r-s)}{(n+r)(p-s)}>1$ implies $8$ cases: $-n<s<p<r$,
$s<-n<r<p$,  $p<r<-n<s$, $r<p<s<-n$, $s<p<r<-n$,  $p<s<-n<r$,
$r<-n<s<p$ and $-n<r<p<s$.

\vskip 3mm \par In \cite{WY}, we proved monotonicity properties of
$\left(\frac{as_r(K)}{n|K|}\right)^{\frac{{n+r}}{r}}$. Here we
prove similar results for mixed $p$-affine surface area.
\begin{proposition}\label{Mono:1}
Let all $K_i\in C^2_+$ be convex bodies with centroid at the
origin.
\vskip 2mm
\noindent
(i) If $-n<r<p$ or $r<p<-n$, one has
\begin{eqnarray*}\label{Mono:2}
\left(\frac{as_p(K_1, \cdots, K_n)}{n\tilde{V}(K_1^\circ, \cdots,
K_n^\circ)}\right)^{n+p} \leq \left(\frac{as_r(K_1, \cdots,
K_n)}{n\tilde{V}(K_1^\circ, \cdots, K_n^\circ)}\right)^{n+r}.
\end{eqnarray*}
\par
\noindent
(ii) If $0<p<r$, or $p<r<-n$, or
$r<-n<0<p$, or $-n<p<r<0$, one has
\begin{eqnarray*}\label{Mono:3}
\left(\frac{as_p(K_1, \cdots, K_n)}{as_0(K_1, \cdots,
K_n)}\right)^{\frac{n+p}{p}} \leq \left(\frac{as_r(K_1, \cdots,
K_n)}{as_0(K_1, \cdots, K_n)}\right)^{\frac{n+r}{r}}.
\end{eqnarray*}
\end{proposition}

\vskip 2mm \noindent{\bf Proof.} \vskip 2mm \noindent(i) We divide
both sides of inequality (\ref{i-3}) by $n\tilde{V}(K_1^\circ,
\cdots, K_n^\circ)$, and get for $\frac{n+p}{n+r}>1$,
\begin{eqnarray}\label{i-4} \frac{as_p(K_1,
\cdots,K_{n})}{n\tilde{V}(K_1^\circ, \cdots, K_{n}^\circ)} \leq
\left( \frac{as_r(K_1, \cdots, K_{n})}{n\tilde{V}(K_1^\circ,
\cdots, K_{n}^\circ)} \right)^{\frac{n+r}{n+p}}.
\end{eqnarray}
Condition $\frac{n+p}{n+r}>1$ implies that $-n< r <p$ or $p< r <-n$.
If $-n<r<p$,  $n+p>0$ and therefore, inequality
(\ref{i-4}) implies inequality (i). If
$p< r<-n$,  $n+p<0$ and therefore,
\begin{eqnarray*}
\left(\frac{as_p(K_1, \cdots, K_n)}{n\tilde{V}(K_1^\circ, \cdots,
K_n^\circ)}\right)^{n+p} \geq \left(\frac{as_r(K_1, \cdots,
K_n)}{n\tilde{V}(K_1^\circ, \cdots, K_n^\circ)}\right)^{n+r}.
\end{eqnarray*}
Switching $r$ and $p$, one obtains the  inequality in
(i):
for $r<p<-n$, $$\left(\frac{as_p(K_1, \cdots,
K_n)}{n\tilde{V}(K_1^\circ, \cdots, K_n^\circ)}\right)^{n+p} \leq
\left(\frac{as_r(K_1, \cdots, K_n)}{n\tilde{V}(K_1^\circ, \cdots,
K_n^\circ)}\right)^{n+r}.$$
\vskip 2mm
\noindent
(ii) Let $s=0$ in inequality (i) of Theorem \ref{ineq2:pmixed}.
Then for $\frac{r(n+p)}{p(n+r)}>1$, $$ as_p(K_1,
\cdots,K_{n}) \leq \big(as_r(K_1,
\cdots,K_{n})\big)^{\frac{p(n+r)}{r(n+p)}} \big(as_0(K_1,
\cdots,K_{n})\big)^{\frac{(r-p)n}{r(n+p)}}.$$
We divide both sides of the
inequality by $as_0(K_1, \cdots, K_n)$ and get $$ \frac{as_p(K_1,
\cdots,K_{n})}{as_0(K_1, \cdots,K_{n})} \leq \left(\frac{as_r(K_1,
\cdots,K_{n})}{as_0(K_1,
\cdots,K_{n})}\right)^{\frac{p(n+r)}{r(n+p)}} .$$
The condition $\frac{r(n+p)}{p(n+r)}>1$
implies that $0<p<r$, or $p<r<-n$, or $-n<r<p<0$, or $r<-n<0<p$.
In the cases $0<p<r$, or $p<r<-n$, or $r<-n<0<p$, one has
$\frac{n+p}{p}>0$ and therefore inequality (ii) holds
true. On the other hand, if $-n<r<p<0$, then $\frac{n+p}{p}<0$ and
hence, $$ \left(\frac{as_r(K_1, \cdots, K_n)}{as_0(K_1, \cdots,
K_n)}\right)^{\frac{n+r}{r}}\leq \left(\frac{as_p(K_1, \cdots,
K_n)}{as_0(K_1, \cdots, K_n)}\right)^{\frac{n+p}{p}}.$$ Switching
$r$ and $p$, one gets inequality (ii):
 if $-n<p<r<0$, then
$$
\left(\frac{as_p(K_1, \cdots, K_n)}{as_0(K_1, \cdots,
K_n)}\right)^{\frac{n+p}{p}}\leq \left(\frac{as_r(K_1, \cdots,
K_n)}{as_0(K_1, \cdots, K_n)}\right)^{\frac{n+r}{r}}.$$

\vskip 3mm
\noindent
Now we treat the case $p=-n$.
The mixed $(-n)$-affine surface area of
$K_1, \cdots, K_n$ is defined as
\begin{equation}\label{def:pmixed:-n} as_{-n}(K_1, \cdots,K_{n})=
\max_{u\in S^{n-1}}\big[f_{K_1}(u) ^{\frac{1}{2n}}
h_{K_1}(u)^{\frac{n+1}{2n}} \cdots f_{K_n}(u)^{\frac{1}{2n}}
h_{K_n}(u)^{\frac{n+1}{2n}}\big].
\end{equation} It is easy to verify that $as_{-n}(K, \cdots,K)$ equals to
$as_{-n}(K)$, the $L_{-n}$ affine surface area of $K$ \cite{MW2}.
We have the following proposition. \vskip 2mm \noindent
\begin{proposition}\label{prop:-n} Let all $K_i$ be convex bodies in $C^2_+$ with centroid at the origin.
Let $p\neq -n$ and $s\neq -n$ be real numbers.
\par\noindent (i) Let $T: \bbR^n \rightarrow\bbR^n $ be an invertible linear
transform. Then
$$as_{-n}(TK_1, \cdots, TK_{n})=|{det}(T)| \
as_{-n}(K_1, \cdots,K_{n}).$$ In particular, if $|det(T)|=1$, then
$as_{-n}(K_1, \cdots, K_n)$ is  affine invariant:  $$
as_{-n}(TK_1, \cdots, TK_{n})= as_{-n}(K_1, \cdots,K_{n}).$$

\vskip 2mm \noindent (ii) Alexandrov-Fenchel type inequalities
$$
as_{-n}^m (K_1, \cdots,K_{n})\leq \prod_{i=0}^ {m-1}as_{-n}(K_1,
\cdots,K_{n-m}, \underbrace{K_{n-i}, \cdots, K_{n-i}}_m),$$ with
equality if the $K_j$,  for $j=n-m+1, \cdots, n$ are dilates.
\newline
In particular, if $m=n$,
\begin{eqnarray*}\label{Holder:pmixed:-n}
as_{-n}^n (K_1, \cdots, K_{n}) \leq  as_{-n}(K_1) \cdots
as_{-n}(K_n).
\end{eqnarray*}
\par
\noindent (iii) If $\frac{n(s-p)}{(n+p)(n+s)}\geq 0$, then
\begin{equation*} as_p(K_1, \cdots, K_n)\leq
\big(as_{-n}(K_1, \cdots, K_n)\big)^{\frac{2n(s-p)}{(n+p)(n+s)}} \
as_s(K_1, \cdots, K_n).
\end{equation*}
\par
\noindent (iv) If $\frac{n(s-p)}{(n+p)(n+s)}\leq 0$, then
\begin{equation*} as_p(K_1, \cdots, K_n)\geq
\big(as_{-n}(K_1, \cdots, K_n)\big)^{\frac{2n(s-p)}{(n+p)(n+s)}} \
as_s(K_1, \cdots, K_n).
\end{equation*}
\end{proposition}

\vskip 2mm
\noindent
{\bf Proof.}
\par
\noindent (i)  By formula (\ref{Curvature:K:TK:2}),
$f_{-n}(TK,v)={det}(T)^2\ f_{-n}(K,u)$. Therefore
\begin{eqnarray*}as_{-n}(TK_1, \cdots, TK_n)&=&\mbox{max}_{v\in
S^{n-1}}[f_{-n}(TK_1,v)\cdots
f_{-n}(TK_n,v)]^{\frac{1}{2n}}\\
&=&  |{det}(T)| \  \mbox{max}_{u\in S^{n-1}}[f_{-n}(K_1,u)\cdots
f_{-n}(K_n,u)]^{\frac{1}{2n}}\\ &=& |{det} (T)|\ \ as_{-n}(K_1,
\cdots, K_n).\end{eqnarray*}
\par
\noindent
(ii) Let $\tilde{g}_0(u)=f_{K_1}(u) ^{\frac{1}{2n}}
h_{K_1}(u)^{\frac{n+1}{2n}} \cdots f_{K_{n-m}}(u)^{\frac{1}{2n}}
h_{K_{n-m}}(u)^{\frac{n+1}{2n}}$ and
$\tilde{g}_{i+1}(u)=f_{K_{n-i}}(u) ^{\frac{1}{2n}}
h_{K_{n-i}}(u)^{\frac{n+1}{2n}}$,  for $i=0, \cdots, m-1$. Then
\begin{eqnarray*}
as_{-n}(K_1, \cdots, K_n)&=&\max_{u\in S^{n-1}} \tilde{g}_0(u)
\tilde{g}_1(u)\cdots \tilde{g}_{m}(u)\\ &\leq & \prod _{i=0}^{m-1}
\left( \max_{u\in S^{n-1}} \tilde{g}_0(u)
\tilde{g}_{i+1}^m(u)\right)^{\frac{1}{m}}\\ &=& \prod_{i=0}^
{m-1}as_{-n}^{\frac{1}{m}}(K_1, \cdots,K_{n-m},
\underbrace{K_{n-i} \cdots, K_{n-i}}_m).
\end{eqnarray*}
Equality holds if and only if for all $i$, $0\leq i \leq  m-1$,
$\tilde{g}_0(u)\tilde{g}_{i+1}^m(u)$ attain their maximum at the
same direction $u_0$. This condition holds true if the $K_j$,  for
$j=n-m+1, \cdots, n$,  are dilates.

\vskip 2mm \par \noindent (iii) and (iv) \begin{eqnarray*}
as_p(K_1, \cdots, K_n)&=&\int _{S^{n-1}} [f_p(K_1,u)\cdots
f_p(K_n,u)]^{\frac{1}{n+p}}\,d\s(u) \\ &=& \int _{S^{n-1}}
[f_s(K_1,u)\cdots f_s(K_n,u)]^{\frac{1}{n+s}} \\
& & \ \
\left(h_{K_1}^{\frac{n+1}{2n}}(u)f_{K_1}^{\frac{1}{2n}}(u)\cdots
h_{K_n}^{\frac{n+1}{2n}}(u)f_{K_n}^{\frac{1}{2n}}(u)
\right)^{\frac{2n(s-p)}{(n+p)(n+s)}}\,d\s(u)
\end{eqnarray*} which is  \begin{equation*} \leq \big(as_{-n}(K_1, \cdots, K_n)\big)^{\frac{2n(s-p)}{(n+p)(n+s)}} \
as_s(K_1, \cdots, K_n), \quad if ~\frac{n(s-p)}{(n+p)(n+s)}\geq 0,
\end{equation*}
and \begin{equation*} \geq \big(as_{-n}(K_1, \cdots,
K_n)\big)^{\frac{2n(s-p)}{(n+p)(n+s)}} \ as_s(K_1, \cdots, K_n),
\quad if ~\frac{n(s-p)}{(n+p)(n+s)}\leq 0.
\end{equation*}
\vskip 5mm

\subsection{$i$-th mixed $p$-affine surface
area and related inequalities}

\par\noindent For all $p\geq 1$ and all real $i$, the
$i$-th mixed $p$-affine surface area of $K, L \in C^2_+$
is defined as \cite{Lut1987, Wa2007}
$$
as_{p,i}(K,L)=\int
_{S^{n-1}}f_p(K,u)^{\frac{n-i}{n+p}}f_p(L,u)^{\frac{i}{n+p}}\,d\s(u).
$$
Recall that  $f_p(K,u)=f_K(u)h_K^{1-p}(u)$. Here we further
generalize this definition to all $p\neq -n$ and all $i$. An
analogous definition for the $i$-th mixed $(-n)$-affine surface
area of $K$ and $L$ is
$$
as_{-n,i}(K,L)=\max_{u\in S^{n-1}}\big[f_{K}(u)
^{\frac{n-i}{2n}} h_{K}(u)^{\frac{(n+1)(n-i)}{2n}}
f_{L}(u)^{\frac{i}{2n}} h_{L}(u)^{\frac{(n+1)i}{2n}}\big].
$$
When $i\in \bbN$, $0 \leq i\leq n$, then, for all $p$, the $i$-th mixed $p$-affine
surface area of $K$ and $L$  is
$$
as_{p,i}(K,L)=as_p(\underbrace{K, \cdots, K}_{n-i}, \underbrace{L,
\cdots, L}_i).
$$
Clearly, for all $p$,  $as_{p,0}(K,L)=as_p(K)$, and
$as_{p,n}(K,L)=as_p(L)$.   When $L=B^n_2$, we write
$as_{p,i}(K)$  for $as_{p,i}(K, B^n_2)$. Thus
\begin{eqnarray*}
as_{p,i}(K)&=&\int_{S^{n-1}}f_p(K,u)^{\frac{n-i}{n+p}}\,d\s(u),
\quad \mbox{for $p\neq -n$},\\
as_{p,i}(K)&=&\max_{u\in S^{n-1}}\big[f_{K}(u) ^{\frac{n-i}{2n}}
h_{K}(u)^{\frac{(n+1)(n-i)}{2n}}\big], \quad \mbox{for $p=-n$.}
\end{eqnarray*}
In particular,  $as_{1,-1}(K)=\int _{S^{n-1}} f_K(u)\,d\s(u)$
is the surface area of $K$.

\vskip 3mm
The next proposition and its proof  is similar to Proposition \ref{ineq2:pmixed} and its proof.
Therefore we omit it.
\begin{proposition}
Let $K$ and $L$ be convex bodies in $C^2_+$ with centroid at the
origin. Let  $i\in \bbR$ and $s\neq -n$, $r \neq -n$,  and $p \neq
-n$ be real numbers.
 \vskip 2mm
\noindent
(i) If
$\frac{(n+p)(r-s)}{(n+r)(p-s)}>1$, then
\begin{equation*} as_{p,i}(K,L) \leq \big(as_{r,i}(K,L)\big)^{\frac{(p-s)(n+r)}{(r-s)(n+p)}}
\big(as_{s,i}(K, L)\big)^{\frac{(r-p)(n+s)}{(r-s)(n+p)}}.
\end{equation*}

\par \noindent (ii) If $\frac{n+p}{n+r} >1$,  then
\begin{eqnarray*} as_{p,i} (K, L) \leq \left(as_{r,i}(K,L)\right)^{\frac{n+r}{n+p}} {\left(n\tilde{V}_i(K^\circ,
L^\circ)\right)^{\frac{p-r}{n+p}}} \end{eqnarray*} where
$\tilde{V}_i(K^\circ, L^\circ)=\frac{1}{n}\int _{S^{n-1}}
\frac{1}{h_K(u)^{n-i}\ h_L(u)^i}\,d\s (u)$ for all $i$.

\par
\noindent (iii) If $\frac{n(s-p)}{(n+p)(n+s)}\geq 0$, then
\begin{equation*} as_{p,i}(K, L)\leq
\big(as_{-n,i}(K, L)\big)^{\frac{2n(s-p)}{(n+p)(n+s)}} \
as_{s,i}(K, L).
\end{equation*}
\par
\noindent (iv) If $\frac{n(s-p)}{(n+p)(n+s)}\leq 0$, then
\begin{equation*} as_{p,i}(K, L) \geq
\big(as_{-n,i}(K, L)\big)^{\frac{2n(s-p)}{(n+p)(n+s)}} \
as_{s,i}(K, L).
\end{equation*}

\end{proposition}

\vskip 3mm
\noindent
The following proposition was proved in
\cite{Lut1987, Wa2007} for $p\geq 1$.

\begin{proposition}\label{Holder:i:pmixed1}
Let $K$ and $ L$  be convex
bodies in $C^2_+$  with centroid at the origin. If $j<i<k$ or $k<i<j$ (equivalently,
$\frac{k-j}{k-i}>1$ ), then for all $p$,
$$as_{p,i}(K,L)\leq
as_{p,j}(K,L)^{\frac{k-i}{k-j}}as_{p,k}(K,L)^{\frac{i-j}{k-j}},$$
with equality if
$K$ and $L$ are dilates of each
other.
\newline
In particular,
$$as_{p,i}(K)\leq
as_{p,j}(K)^{\frac{k-i}{k-j}}as_{p,k}(K)^{\frac{i-j}{k-j}},$$ with
equality if $K$ is a ball.

\end{proposition}

\vskip 2mm
\noindent  For $p\neq -n$, the proof is the same as
the proof in \cite{Lut1987, Wa2007}. For $p=-n$, it is similar to
the proof of Proposition \ref{ineq1:pmixed}.  Note that for $i\in
\bbN$, $0<i<m$, $m=j$ and $k=0$, the proposition is a direct
consequence of Proposition \ref{ineq1:pmixed}.

\vskip 2mm \noindent In Proposition \ref{Holder:i:pmixed1}, if
$j=0$ and  $k=n$, then for all $p$ and $0\leq i \leq n$
\begin{eqnarray}\label{Holder:i:pmixed3}as_{p,i}^n (K,L) \leq
as_{p}^{{n-i}}(K)as_{p}^{{i}}(L).\end{eqnarray}

\vskip 2mm \noindent If we let $i=0$ and  $j=n$, then for all
$k\leq 0$ and for all $p$
\begin{eqnarray}\label{Holder:i:pmixed5}as_{p,k}^n(K, L)\geq
as_p^{n-k} (K) as_p^k(L). \end{eqnarray}

\vskip 2mm \noindent Let $i=n$, $j=0$ and $k>n$. Then inequality
(\ref{Holder:i:pmixed5}) also holds true for $k\geq n$ and all
$p$. In both (\ref{Holder:i:pmixed3}) and
(\ref{Holder:i:pmixed5}), equality holds for all $p$ if $K$ and
$L$ are dilates. \vskip 2mm \noindent From inequality
(\ref{Holder:i:pmixed3}) and Corollary 4.1 in \cite{WY}, one gets
that
\begin{eqnarray}\label{Holder:i:pmixed4}as_{p,i}^n(K,L)as_{p,i}^n(K^\circ,L^\circ)&\leq&
\big(as_{p}(K)as_{p}(K^\circ)\big)^{{n-i}} \big(
as_{p}(L)as_{p}(L^\circ)\big)^{{i}}\nonumber\\ &\leq&
n^{2n}(|K||K^\circ|)^{n-i} (|L||L^\circ|)^i\end{eqnarray} holds
true for all $p\geq 0$ and $0< i < n$. The inequality also holds
if $i=0$ and $i=n$ \cite{WY}. We apply  Blaschke-Santal\'{o}
inequality to inequality (\ref{Holder:i:pmixed4}) and get
\begin{eqnarray*}as_{p,i}(K,L)as_{p,i}(K^\circ,L^\circ)\leq
as_p^2(B^n_2)\end{eqnarray*} for all $p\geq 0$ and $0\leq i\leq
n$. Equality holds true if $K$ and $L$ are ellipsoids that are
dilates of each other. Hence we have proved the following
proposition, which, for $p\geq 1$, was proved in \cite{Wa2007}.

\begin{proposition}\label{Holder:i:pmixed7}
Let $K$ and $L$ be convex bodies in $C^2_+$ with centroid at the origin.
If $p\geq 0$ and $0\leq i\leq n$, then
$$as_{p,i}(K,L)as_{p,i}(K^\circ,L^\circ)\leq as_p^2(B^n_2),$$
with equality if $K$ and $L$ are ellipsoids that are dilates of
each other.
\end{proposition}

\vskip 2mm We now establish isoperimetric inequalities for
$as_{p,i}(K)$.

\begin{proposition}\label{Holder:i:pmixed6} Let $K\in C^2_+$ be a
convex body with centroid at the origin.

\par\noindent (i) If $p\geq 0$ and $0\leq i\leq n$, then
$$\frac{as_{p,i}(K)}{as_{p,i}(B^n_2)}\leq
\left(\frac{|K|}{|B_2^n|}\right)^{\frac{(n-p)(n-i)}{(n+p)n}}
$$ with equality if $K$ is a ball. Moreover,
$as_{p,i}(K)as_{p,i}(K^\circ)\leq as_p^2(B^n_2)$ with equality if
$K$ is a ball.

\par\noindent (ii) If $p\geq 0$ and $i\geq n$, then
$$\frac{as_{p,i}(K)}{as_{p,i}(B^n_2)}\geq
\left(\frac{|K|}{|B_2^n|}\right)^{\frac{(n-p)(n-i)}{(n+p)n}}
$$ with equality if $K$ is a ball. Moreover,
$as_{p,i}(K)as_{p,i}(K^\circ)\geq as_p^2(B^n_2)$ with equality if
$K$ is a ball.

\par\noindent (iii) If $-n<p<0$ and $i\leq 0$, then
$$\frac{as_{p,i}(K)}{as_{p,i}(B^n_2)}\geq
\left(\frac{|K|}{|B_2^n|}\right)^{\frac{(n-p)(n-i)}{(n+p)n}}
$$ with equality if $K$ is a ball. Moreover,
$as_{p,i}(K)as_{p,i}(K^\circ)\geq c^{n-i} as_p^2(B^n_2)$ where $c$
is the universal constant in the inverse Santal\'{o} inequality
\cite{BM, GK2}.

\par\noindent (iv) If $p<-n$ and $i\leq 0$, then
$$\frac{as_{p,i}(K)}{as_{p,i}(B^n_2)}\geq c^{\frac{(n-i)p}{n+p}}
\left(\frac{|K|}{|B_2^n|}\right)^{\frac{(n-p)(n-i)}{(n+p)n}}.
$$ Moreover,
$as_{p,i}(K)as_{p,i}(K^\circ)\geq c^{n-i} as_p^2(B^n_2)$ where $c$
is the same constant as in (iii).

\par\noindent (v) If $i\leq 0$, then
$$\frac{as_{-n,i}(K)}{as_{-n,i}(B^n_2)}\geq
\left(\frac{|K|}{|B^n_2|}\right)^{\frac{n-i}{n}}.$$ Moreover,
$as_{-n,i}(K)as_{-n,i}(K^\circ)\geq as_{-n,i}^2(B^n_2).$
\end{proposition}

\par \noindent
{\bf Proof.}

\vskip 2mm \par \noindent (i) For  $i=n$, the equality holds
trivially. For $i=0$, the inequality was proved in \cite{WY}. We
now prove the case $0<i<n$. $L=B^n_2$ in inequality
(\ref{Holder:i:pmixed3}) gives
\begin{equation}\label{Holder:i:pmixed2}\left(\frac{as_{p,i}(K)}{as_{p,i}(B^n_2)}\right)^n\leq
\left(\frac{as_{p}(K)}{as_p(B^n_2)}\right)^{{n-i}}\end{equation}
for all $p\neq -n$ and $0\leq i \leq n$. We also use that
$as_{p,i}(B^n_2)=as_p(B^n_2)$.  Then,  as $as_p(B^n_2)=n|B^n_2|$,
we get
 for all $p\geq 0$ and $0\leq i\leq n$,
the following isoperimetric inequality  as
a consequence of the $L_p$ affine isoperimetric inequality in
\cite{WY}
$$\frac{as_{p,i}(K)}{as_{p,i}(B^n_2)}\leq
\left(\frac{as_{p}(K)}{as_p(B^n_2)}\right)^{\frac{{n-i}}{n}}\leq
\left(\frac{|K|}{|B_2^n|}\right)^{\frac{(n-p)(n-i)}{(n+p)n}},
$$ with equality if $K$ is a ball. The inequality
${as_{p,i}(K)as_{p,i}(K^\circ)}\leq {as_{p,i}^2(B^n_2)}$ follows
from Proposition \ref{Holder:i:pmixed7} with $L=B^n_2$.

\vskip 2mm \par \noindent (ii) For  $i=n$, the equality holds
trivially. Similarly, let $L=B^n_2$ in inequality
(\ref{Holder:i:pmixed5}), then for all $p\neq -n$, and $i\geq n$
or $i\leq 0$,
\begin{equation}\label{Holder:i:pmixed8}\left(\frac{as_{p,i}(K)}{as_{p,i}(B^n_2)}\right)^n\geq
\left(\frac{as_{p}(K)}{as_p(B^n_2)}\right)^{{n-i}}.\end{equation}

\par\noindent Hence for $i\geq n$ and $p\geq 0$, the $L_p$
affine isoperimetric inequality in \cite{WY} implies that
$$\frac{as_{p,i}(K)}{as_{p,i}(B^n_2)}\geq
\left(\frac{as_{p}(K)}{as_p(B^n_2)}\right)^{\frac{{n-i}}{n}}\geq
\left(\frac{|K|}{|B_2^n|}\right)^{\frac{(n-p)(n-i)}{(n+p)n}}
$$ with equality if $K$ is a ball. Moreover, by
Corollary 4.1 (i) in \cite{WY} and the remark after it, one has
for all $i \geq n$
 $$\left(\frac{as_{p,i}(K)as_{p,i}(K^\circ)}{as_{p,i}^2(B^n_2)}\right)^n\geq
\left(\frac{as_{p}(K)as_{p}(K^\circ)}{as_p^2(B^n_2)}\right)^{{n-i}}\geq
1,$$  or equivalently, $as_{p,i}(K)as_{p,i}(K^\circ)\geq
{as_{p,i}^2(B^n_2)}$, with equality if $K$ is a ball.

\vskip 2mm \noindent (iii) If $i\leq 0$ and $-n< p<0$, inequality
(\ref{Holder:i:pmixed8}) and Theorem 4.2 (ii) of \cite{WY} imply
that $$\frac{as_{p,i}(K)}{as_{p,i}(B^n_2)}\geq
\left(\frac{as_{p}(K)}{as_p(B^n_2)}\right)^{\frac{{n-i}}{n}}\geq
\left(\frac{|K|}{|B_2^n|}\right)^{\frac{(n-p)(n-i)}{(n+p)n}}
$$ with equality if $K$ is a ball.  By Corollary 4.1 (ii) of \cite{WY} and the remark after it,
$$\left(\frac{as_{p,i}(K)as_{p,i}(K^\circ)}{as_{p,i}^2(B^n_2)}\right)^n\geq
\left(\frac{as_{p}(K)as_{p}(K^\circ)}{as_p^2(B^n_2)}\right)^{{n-i}}\geq
c^{n(n-i)}, $$ or equivalently,  $as_{p,i}(K)as_{p,i}(K^\circ)\geq
 c^{n-i}{as_{p,i}^2(B^n_2)}$ where $c$ is the constant in the inverse Santal\'{o} inequality \cite{BM,GK2}.

\vskip 2mm \noindent (iv) If $i \leq 0$ and $p<-n$  inequality
(\ref{Holder:i:pmixed8}) and Theorem 4.2 (iii) of \cite{WY} imply that
$$\frac{as_{p,i}(K)}{as_{p,i}(B^n_2)}\geq
\left(\frac{as_{p}(K)}{as_p(B^n_2)}\right)^{\frac{{n-i}}{n}}\geq
c^{\frac{(n-i)p}{n+p}}
\left(\frac{|K|}{|B_2^n|}\right)^{\frac{(n-p)(n-i)}{(n+p)n}}.
$$  The proof of  $as_{p,i}(K)as_{p,i}(K^\circ)\geq
 c^{n-i}{as_{p,i}^2(B^n_2)}$ is same as in (iii).

\vskip 2mm\noindent (v) Inequality (\ref{Holder:i:pmixed5})
implies that $as_{-n,i}(K)^n\geq as_{-n}(K)^{n-i}$ for $i\leq 0$. As
$as_{-n,i}(B^n_2)=1$ for all $i$,
 $$\left(\frac{as_{-n,i}(K)}{as_{-n,i}(B^n_2)}\right)^n \geq \left(\frac{as_{-n}(K)}{as_{-n}(B^n_2)}\right)^{n-i}
 \geq \left(\frac{|K|}{|B^n_2|}\right)^{n-i}.$$ The second
 inequality follows from the $L_{-n}$ affine isoperimetric
 inequality in \cite{WY}.

\vskip 2mm  \par \noindent Moreover, by Corollary 4.2 in
\cite{WY},
 for $i\leq 0$,
$$\left(\frac{as_{-n,i}(K)as_{-n,i}(K^\circ)}{as_{-n,i}^2(B^n_2)}\right)^n
 \geq \left(\frac{as_{-n}(K)as_{-n}(K^\circ)}{as_{-n}^2(B^n_2)}\right)^{n-i}
 \geq 1, $$ or equivalently, $as_{-n,i}(K)as_{-n,i}(K^\circ)\geq
 as_{-n,i}^2(B^n_2)$.

\vskip 2mm \par \noindent {\bf Remark.} The example $K(R,
\varepsilon)$ mentioned in the remarks after Proposition
\ref{Iso:pmixed:2} shows that we  cannot expect to get
strictly positive lower bounds in (i) of Proposition
\ref{Holder:i:pmixed6} for $p>0$ and $0\leq i< n$. In fact, by
inequality (\ref{Holder:i:pmixed2}), one has
$$\left(\frac{as_{p,i}(K(R,
\varepsilon))}{as_{p,i}(B^n_2)}\right)^n\leq
\left(\frac{as_{p}(K(R,
\varepsilon))}{as_p(B^n_2)}\right)^{{n-i}}.$$ As in \cite{WY}, $as_{p}(K(R, \varepsilon))\rightarrow 0$ for $p>0$ as
$R\rightarrow \infty$ and $\varepsilon \rightarrow 0$. $0\leq i
<n$ implies that $n-i>0$, and therefore $as_{p,i}(K(R,
\varepsilon)) \rightarrow 0$.

\vskip 2mm \par This example also shows that, likewise, we
cannot expect finite upper bounds in $(ii)$ (for $p>0$ and
$i>n$), $(iii)$ (for $-n<p<0$ and $i\leq 0$), and $(iv)$ (for
$p<-n$ and $i\leq 0$), of Proposition \ref{Holder:i:pmixed6}. For
instance, if $i\leq 0$, by inequality (\ref{Holder:i:pmixed8}),
one has
$$\left(\frac{as_{p,i}(K(R, \varepsilon))}{as_{p,i}(B^n_2)}\right)^n\geq
\left(\frac{as_{p}(K(R,\varepsilon))}{as_p(B^n_2)}\right)^{{n-i}}.$$
For $-2<p<0$, one has $as_{p}(K(R,\varepsilon))\rightarrow \infty$
as $R\rightarrow \infty$ and $\varepsilon \rightarrow 0$.
Therefore, if $i\leq 0$, we obtain that $as_{p,i}(K(R,
\varepsilon))\rightarrow \infty$ as $R\rightarrow \infty$ and
$\varepsilon \rightarrow 0$, i.e., there are no finite upper
bounds in $(iii)$.

\vskip 2mm \par \noindent {\bf Remark.} In $(iv)$, if $p=-\infty$,
then for all $i\leq 0$,  $$\int \int_{S^{n-1}\times S^{n-1}}
\big(h_K(u)h_{K^\circ}(v)\big)^{i-n}\,d\s (u) \,d\s (v)\geq
 c^{n-i}{as_{p,i}^2(B^n_2)}$$ or equivalently, for all $i\leq 0$,  $$\int \int_{S^{n-1}\times S^{n-1}}
\big(\r _K(u)\r _{K^\circ}(v)\big)^{n-i}\,d\s (u) \,d\s (v)\geq
 c^{n-i}{as_{p,i}^2(B^n_2)}.$$ In particular, if $i=0$, this is equivalent to the
 inverse Santal\'{o} inequality \cite{BM}.

\section{Illumination surface bodies}
We now  define a new  family of bodies associated with  a given
convex body $K$. These new bodies are a variant of the
illumination bodies \cite{W1} (compare also \cite{ SW5}).

\bd  ({\bf Illumination surface body}) Let $s\geq 0$ and
$f:\partial K\rightarrow \bbR$ be a nonnegative, integrable
function. The illumination surface body $K^{f,s}$ is defined as
$$\displaystyle K^{f,s}=\left\{x: \mu _f({\partial K \cap
\overline{[x,K]\backslash K}})\leq s\right\}.$$ \ed

\vskip 2mm \par Obviously, $K\subseteq K^{f,s}$ for any $s\geq 0$
and any nonnegative, integrable function $f$.  Moreover,
$K^{f,s}\subseteq K^{f,t}$ for any $0\leq s\leq t$. \vskip 2mm
Notice also that $K^{f,s}$ needs to be neither bounded nor convex:

\begin{example}\label{ex-1}
\par \noindent
Let $K=B_{\infty}^2=\{x \in \mathbb{R}^2: \mbox{max}_{1 \leq i
\leq 2} |x_i| \leq 1 \}$ and
\begin{equation*}
f(x)=\left\{
\begin{array}{cc}
\frac{1}{12}, & x\in [(-1,1), (1,1)]\cup [(1,1), (1,-1)]\\
\frac{1}{6}, & ~\mbox{otherwise}~
\end{array} \right.
\end{equation*}
\par \noindent $K^{f,s}$ is calculated as follows. If $s<\frac{1}{6}$, then $K^{f,s}=K$.
\par \noindent
If $s\in [\frac{1}{6}, \frac{1}{3})$, then $$K^{f,s}=\{(x_1,x_2):
x_1\geq -1, x_2\in [-1,1]; \ \mbox{or} \  x_2\geq -1, x_1\in
[-1,1]\}.$$
\par \noindent If $s\in [\frac{1}{3}, \frac{1}{2})$, then
$$K^{f,s}=\{(x_1,x_2): x_1, x_2\geq -1;\ \mbox{or} \ x_1\leq -1, x_2\in [-1,1];
 \ \mbox{or} \ x_2\leq -1, x_1\in [-1,1] \}.$$
If $s\in[\frac{1}{2}, \frac{2}{3})$, then $K^{f,s}=\{(x_1,x_2):
x_1\geq -1 ~or~ x_2 \geq -1\}.$
\par \noindent Except for $s<1/6$, all of
them are neither bounded nor convex.
\par \noindent If $s\geq \frac{2}{3}$, then $K^{f,s}=\bbR^2$
\end{example}
\vskip 2mm \noindent The following lemmas describe some of the
properties of the bodies $K^{f,s}$.

\bl \label{Le:1} Let $s\geq 0$ and $f : \partial K \rightarrow
\mathbb{R}$ be a nonnegative, integrable function. Then
\par
\noindent
(i)\ \ $K^{f,s}=\bigcap _{\delta
> 0} K^{f,s+\delta}$.
\par
\noindent
(ii)\ \  $K^{f,s}$ is star convex, i.e.,  for all $x \in K^{f,s}$:
$[0,x]  \subset   K^{f,s}$. \el
\par
\noindent {\bf Proof.}
\newline
(i) We only need to show that
$K^{f,s}\supseteq \bigcap _{\delta
> 0} K^{f,s+\delta}$. Let  $x\in \bigcap _{\delta
> 0} K^{f,s+\delta}$. Then for all $\delta >0$,
$ \mu_f({\partial K \cap \overline{ [x,K]\backslash K}}) \leq
s+\delta$. Thus, letting  $\delta \rightarrow 0$, $
\m_f({\partial K \cap \overline{ [x,K]\backslash K}})\leq s$.

\par\noindent (ii)  Let  $x\in K^{f,s}$. We claim that
$[0,x]\subseteq K^{f,s}$. Let  $y\in [0,x]$. Since $y\in
[0,x]\subset  [x,K]$, we have $ [y,K]\setminus K\subseteq
[x,K]\setminus K$ and thus $\partial K\cap \overline{
[y,K]\setminus K}\subseteq \KC.$ This implies that
$$ \mu_f({\partial K\cap \overline{ [y,K]\setminus K}})\leq
\mu _f({\partial K\cap \overline{ [x,K]\setminus K} })\leq s
$$ and hence  $y\in K^{f,s}$.

\vskip 2mm \noindent {\bf Remark.} We {\em can not} expect
$K^{f,s}$ to be convex,  even for $K=\bn n$ and $f$  smooth.
Indeed, let $K=B_2^2$ and $s=\frac{1}{64}$. Define
$$
f(x)=\left\{\begin{array}{ll} \frac{1}{4\pi}, & \mbox{$x$ is in the first and third quadrant}\\
\frac{1}{16 \pi}, & \mbox{$x$ is in the fourth quadrant}\\
\frac{23}{16\pi}, & \mbox{$x$ is in the second
quadrant}\end{array} \right.
$$ Then $K^{f,s}$ is not convex. In
fact, $K^{f,s}$ contains the arc from the point $(
\tan(\frac{\pi}{32}), 1)$ to the point $(1, \tan(\frac{\pi}{32}))$
of the Euclidean ball centered at $0$ with radius
$r=\mbox{sec}(\frac{\pi}{32})$. Moreover, the point
$(\mbox{sec}(\frac{\pi}{20}), 0)$ is on the boundary  of
$K^{f,s}$. The tangent line at $(1, \tan(\frac{\pi}{32}))$ of
$B^2_2(0,r)$ is $y=\frac{-x +r^2}{\sqrt{r^2-1}}$. This tangent
line intersects the  $x$-axis at $(r^2, 0)=
(\mbox{sec}^2(\frac{\pi}{32}), 0)$. Since
$\mbox{sec}^2(\frac{\pi}{32}) \sim 1.009701<
\mbox{sec}(\frac{\pi}{20}) \sim1.01264$, $K^{f,s}$ is not convex.
\par
We can modify  $f$ so that it becomes smooth also at the points
$(\pm 1, 0)$ and $(0, \pm 1 )$  and $\partial K^{f,s}$ still
intersects the positive $x$-axis at the point
$(\mbox{sec}(\frac{\pi}{20}), 0)$. Therefore, $K^{f,s}$ is not
convex, even if  $f$ is smooth. \vskip 3mm

\bl \label{Le:2}
Let $s\geq 0$ and $f: \partial K \rightarrow \mathbb{R}$ be an  integrable, $\m
_K$-almost everywhere strictly positive function. Then \\
(i) $K=K^{f,0}$ .\\
(ii) There exists $s_0>0$, such that for all $0\leq s\leq s_0$,
$K^{f,s}$ is bounded. \el
\par \noindent
{\bf Proof.} (i)  It is enough to prove that $K^{f,0}\subseteq K$.
Suppose this is not the case.  Then there is  $x\in K^{f,0}$ but
$x\notin K$. Since $0\in \mbox{int}(K)$, there is $\a>0$ such that
\begin{equation}\label{alpha}
B_2^n(0,\a)\subseteq K\subseteq B_2^n(0,1/\a).
\end{equation}
Let $y\in [0,x]\cap \partial K $ and
$\mbox{Con}(x,\a)=[x,B_2^n(0,\a)]$ be the convex hull of $x$ and
$B_2^n(0,\a)$. $H(y,N_K(y))\cap \mbox{Con}(x,\a)$ contains a
$(n-1)$-dimensional Euclidean ball with radius (at least)
$r_1=\a\frac{\|x-y\|}{\|x\|}>0$.

\vskip 3mm \par \noindent Hence $\m_K(\pt K \cap \overline{
[x,K]\setminus K})\geq |H(y,N_K(y))\cap Con(x,\a)| \geq r_1^{n-1}
\  |B_2^{n-1}|\
> 0.$ Let
$$E_j=\left\{ z\in \KC: f(z)\geq \frac{1}{j}\right\}, \ \ j=1,2,\cdots .$$
As $\mu_K \left(\{ z \in \partial K: f(z)=0\}\right)=0$ and $E_j
\subseteq E_{j+1}$ for all $j$,
$$
\mu_K\big(\pt K \cap \overline{ [x,K]\setminus K}\big) =
\mu_K\big( \bigcup_{j=1}^{\infty}E_j \big) = \lim_{j \rightarrow
\infty} \mu_K(E_j).
$$
Therefore there exists  $j_1$ such that $\mu_K(E_{j_1}) >0$. Thus
\begin{eqnarray*}
\m _f({\KC})\geq \m _f({E_{j_1}})\geq \frac{\mu_K
(E_{j_1})}{j_1}>0
\end{eqnarray*}
which  contradicts that  $x\in K^{f,0}$.

\par
\noindent
(ii) is an immediate consequence of Lemma \ref{Le:1} (i) and Lemma
(\ref{Le:2}) (i). Indeed, these lemmas  imply
that $K=K^{f,0}=\bigcap _{s
> 0} K^{f,s}$. So, also using  (\ref{alpha}), there exists $s_0>0$ such that for all
$0\leq s \leq s_0$, $K^{f,s}\subset 2K\subset B^n_2(0,
\frac{2}{\a})$. In particular, $K^{f,s_0}\subset B^n_2(0,
\frac{2}{\a}).$

\vskip 2mm \noindent {\bf Remark.} The assumption that $f$ is
$\m_K$-almost everywhere strictly positive is necessary in order
that $K^{f,0}=K$. To see that,  let $K=B_2^2$ and
$$\displaystyle f(x,y)=\left\{
\begin{array}{cc}
0&~  x=\sqrt{1-y^2}, y\in [-1,1],\\
\frac{1}{\pi}& ~ \mbox{otherwise}.~
\end{array}\right .$$
Then $K^{f,0}=K\cup \{(x,y): x\geq 0, |y|\leq 1\}$. \vskip 2mm
\par \noindent This example also shows that there is no $s_0$ such
that $K^{f,s}$ is bounded for all $0\leq s \leq s_0$ unless $f$ is
$\m_K$-almost everywhere strictly positive.

\vskip 3mm  \par \noindent Let $K$ be a convex body with $0\in
\mbox{int}(K)$. Let $f: \partial K \rightarrow \mathbb{R}$ be an
integrable, $\m _K$-almost everywhere strictly positive function.
For $x\notin K$, let $t_0=t_0(x)$ be the strictly positive real
number such that $t_0x=\pt K \cap [0,x]$. Define $h_x(t)$ to be
$$
h_x(t)=\m_f({\partial K\cap \overline{ [tx,K]\setminus K}}), \quad
\mbox{$t\geq t_0$}.
$$
Clearly $h_x(t_0)=0$. Moreover $h_x(t)\leq s$ if $tx\in K^{f,s}$,
and $h_x(t)>s$ if $tx\notin K^{f,s}$.

\bl  \label{Lemma:prop} Let $K$ be a convex body  in
$\mathbb{R}^n$ and $f: \partial K \rightarrow \mathbb{R}$ be an
integrable, $\m _K$-almost everywhere strictly positive function.
\vskip 1mm \par \noindent (i) $h_x(t)$ is increasing and left
continuous on $[t_0, \infty)$.

\vskip 1mm \par \noindent (ii) $K^{f,s}$ is closed for all $s\geq
0$. In particular, it is compact for all  $0\leq s\leq s_0$.

\vskip 1mm \par \noindent If $K$ is in addition strictly convex,
then

\vskip 1mm \par \noindent (iii) $h_x(t)$ is continuous on
$[t_0,\infty)$.

\vskip 1mm \par \noindent (iv) For any $0 \leq s \leq  s_0$ and
$x\in \partial K^{f,s}$, one has $ \m_f({\partial K\cap \overline{
[x,K]\setminus K}})=s. $ \el

\noindent {\bf Proof.} \vskip 2mm \par \noindent (i) If $t_1\leq
t_2$, then $\partial K\cap \overline{ [t_1 x,K]\setminus
K}\subseteq
\partial K\cap \overline{ [t_2x,K]\setminus K}$. Thus
$h_x(t_1)\leq h_x(t_2)$.

\vskip 2mm \par \noindent Let now $t>t_0$ and $(t_m)_{m \in
\mathbb{N}}$ be a sequence,  increasing to $t$. Then, by
monotonicity of $h_x$, $h_x(t_m) \leq h_x(t)$ for all $m$ and thus
$\mbox{lim}_m h_x(t_m) \leq h_x(t)$.  We have to show that
$\mbox{lim}_m h_x(t_m) \geq h_x(t)$. Let $y\in
\mbox{relint}_{\partial K}(\partial K\cap \overline{
[tx,K]\setminus K})$, where $\mbox{relint}_B(A)$ is the relative
(with respect to $B$) interior of a set $A \subseteq B$, i.e.,
$\mbox{relint}_B(A) =\{x\in A: \mbox{there is a $\d>0$, such that
$B(x,\d)\cap B\subset A$}\}$. Then $y\in \mbox{int} \left( [tx, K]
\right)$, and therefore  there exists  $m_0(y)\in \mathbb{N}$,
such that $y\in \mbox{int}([t_{m_0(y)}x, K])$. This implies that
$y \in \overline{[t_{m_{0}(y)}x, K] \setminus K}\cap \partial K$
and thus
 $$\mbox{relint}_{\partial K} (\overline{[tx, K]\setminus
K}\cap \partial K)\subset \bigcup _{m\geq 1} \overline{[t_{m}x,
K]\setminus K}\cap
\partial K.
$$
By  continuity of the measure $\mu_f$ from below, one has
\begin{eqnarray*}
h_x(t) &=& \mu_f \big(\overline{[tx, K]\setminus K} \cap \partial
K\big) = \mu_f \left(\mbox {relint}_{\partial K}
\big(\overline{[tx, K]
\setminus K} \cap \partial K\big)\right)\\
&\leq& \mu_f\bigg( \bigcup_{m \geq 1} \big( \overline{[t_m x,
K]\setminus K}\cap
\partial K\big)\bigg)\\  &=& \lim_m \mu _f \left(\partial K\cap
\overline{ [t_mx,K]\setminus K}\right) = \mbox{lim}_m h_x(t_m)
\end{eqnarray*}

\vskip 2mm \noindent (ii) It will follow from Lemma \ref {Le:2}
(ii) that $K^{f,s}$ is compact for  $0 \leq s \leq s_0$, once we
have proved that $K^{f,s}$ is closed.
\par To that end, we show that  $(K^{f,s})^c$, the
complement of $K^{f,s}$ in $\mathbb{R}^n$,  is open for all $s\geq
0$. Suppose this is not the case. Then there exists $x\in
(K^{f,s})^c$ and a sequence $(x_m)_{m \in \mathbb N}$, such that
$x_m\rightarrow x$ as $m\rightarrow \infty$ but $x_m\in K^{f,s}$
for all $m$. Without loss of generality, we can assume that $x_m$
are not in the ray of $\{tx: t\geq 0\}$. Otherwise, if $x_m \in
K^{f,s}$ are in the ray, then $h_x
\left(\frac{\|x_m\|}{\|x\|}\right)\leq s$ and by (i), $\lim _m h_x
\left(\frac{\|x_m\|}{\|x\|}\right) =h_x(1)\leq s$. This
contradicts with $h_x(1)>s$.

Now we let $K_m= [x_m, K]$. For  sufficiently big $m$,
$\pt K_m \cap [0,x] \neq \emptyset$. Suppose not, then $x\in K_m$
implies that $[x,K]\subset K_m$, and hence
$\overline{[x,K]\setminus K}\cap \pt K \subset
\overline{K_m\setminus K}\cap \pt K$. Since $\m _f
(\overline{[x,K]\setminus K}\cap \pt K)>s$, one gets
$\m_f(\overline{K_m\setminus K}\cap \pt K)>s$, a contradiction
with $x_m\in K^{f,s}$. Let $y_m=\pt K_m \cap[0,x]$. Thus $\mu_f ( \overline{[y_m, K]\setminus K}
\cap \partial K) \leq s$. Let $\alpha$ be
as in (\ref{alpha}). Similarly, $\partial
\big([x_m, B^n_2(0,\a)]\big) \cap [0,x] \neq \emptyset$ for  sufficiently
big $m$ and we denote $z_m=\partial
\big([x_m, B^n_2(0,\a)]\big) \cap [0,x]$.
\par It is easy to check that $0\leq \|x\|-\|y_m\|\leq
\|x\|-\|z_m\|$ for any $m$. As $\a\leq \|z_m\|\leq \|x\|$ and
$\frac{\a}{\|z_m\|}\leq \frac{\|x-x_m\|}{\|x\|-\|z_m\|}$, one has
$\|x\|-\|z_m\| \leq \frac{\|x\| \ \|x-x_m\|}{\a}$. Thus
$z_m\rightarrow x$, and hence also $y_m\rightarrow x$, as
$m\rightarrow \infty$. Therefore we can choose a subsequence
$(y_{m_k})_{k \in \mathbb N}$ that is monotone increasing to $x$.
 By (i) with $t_{m_k}=\frac{\|y_{m_k}\|}{\|x\|}$,
$h_x(t_{m_k})\nearrow h_x(1)$ as
$k\rightarrow  \infty$.
Since for all $k$, $h_x(t_{m_k})\leq s$,
one has  $h_x(1)\leq s$, a contradiction.
\par
\vskip 2mm \noindent (iii) It is enough to prove that $h_x(t)$ is
right continuous on $[t_0, \infty)$. To do so, let $t \geq t_0$
and let
 $(t_m)_{m \in \mathbb {N}}$ be a sequence decreasing to $t$.
 By (i), $h_x(t_m) \geq h_x(t)$ for all $m$,  thus $\mbox{lim}_m h_x(t_m)
\geq h_x(t)$ and we have to show that  $\mbox{lim}_m h_x(t_m)
\leq h_x(t)$.
 We claim that if $K$ is strictly convex, then
 \begin{equation}\label{rightcont}
\partial K\cap \overline{ [tx,K]\setminus K} = \bigcap
_{m=1}^\infty \big[\partial K\cap \overline{ [t_m x,K]\setminus
K}\big].
\end{equation}
We only need to prove that $\bigcap _{m=1}^\infty \big[\partial
K\cap \overline{ [t_m x,K]\setminus K}\big] \subseteq \partial
K\cap \overline{ [tx,K]\setminus K}$. Let $z_0\in \bigcap
_{m=1}^\infty \big[\partial K\cap \overline{ [t_m x,K]\setminus
K}\big]$. Thus $z_0\in \pt K$. Let $l(z_0,tx)$ be the
line passing through $tx$ and $z_0$. We have two cases.\\
Case 1: $l(z_0,tx)$ is in a  tangent hyperplane of $K$. Then
$l(z_0,tx) \cap \partial K=\{z_0\}$ by strict convexity of $K$.
Therefore, $\{z_0\}=\overline{[z_0, tx]\setminus K}\cap
\partial K \subseteq \overline{[tx, K]\setminus K}\cap \partial
K$.
\\ Case 2: $l(z_0,tx) \cap \partial K$ consists of two points,
$z_0$ and $z_1$. As $z_0\in \bigcap _{m=1}^\infty \big[\partial
K\cap \overline{ [t_m x,K]\setminus K}\big]$, we must have
$\|tx-z_0\|<\|tx-z_1\|$. Therefore, $\{z_0\}=\overline{[z_0,
tx]\setminus K}\cap \partial K \subseteq \overline{[tx,
K]\setminus K}\cap \partial K$.
\par
\noindent
Hence by  (\ref{rightcont})  and continuity of the measure $\mu_f$ from above,
\begin{eqnarray*} \nonumber
h_x(t)&=&\mu_f \bigg(\bigcap _{m=1}^\infty \left[\partial K\cap
\overline{ [t_m x,K]\setminus K}\right]\bigg)\\
&=& \lim_{m}
\mu_f \left(
 \partial K\cap \overline{ [t_m x,K]\setminus K}\right)= \lim_{m} h_x(t_m).
\end{eqnarray*}

\noindent (iv) Let $0\leq s\leq s_0$, and $x\in
\partial K^{f,s}$ which implies that $h_x(1) \leq
s.$ Define $\Phi_x (s)=\{t: h_x(t)=s\}.$ Then $\Phi_x (s) \neq
\emptyset$. Indeed, let $t_a x=\partial
B^n_2(0,\frac{3}{\alpha})\cap T_x$, where $\alpha$ is as in
(\ref{alpha}) and $T_x=\{tx: t\geq t_0(x)>0\}$. The proof of Lemma
\ref{Le:2} (ii) shows that $K^{f,s_0}\subset B^n_2(0,
\frac{2}{\a})$. It is clear that $t_{\a}x\notin K^{f,s_0}$, and
hence $h_x(t_{\a})>s_0$. In fact, if $t_{\a}x\in K^{f,s_0}$, then
$t_{\a}x\in B^n_2(0,\frac{2}{\alpha})$, but by definition of
$t_{\a}$, $t_{\a}x\in
\partial B^n_2(0,\frac{3}{\alpha})$. This is a contradiction.

By continuity of $h_x(\cdot)$, there must exist $t\in [t_0,
t_{\a}]$, such that $h_x(t)=s$. This also shows that $\bar{t}=\sup
\Phi_x (s) \leq t_{\a }$.  Clearly $h_x(\bar{t})=s$ and thus
$\bar{t}x\in K^{f,s}.$ This implies that $\bar{t}\leq 1$ because
$x\in \partial K^{f,s}$. Suppose $\bar{t}<1$. Then
$s=h_x(\bar{t})\leq h_x(1)\leq s$ by monotonicity of $h_x(\cdot)$,
a contradiction with $\bar{t}=\sup \Phi (s)$. Thus $\bar{t}=1$ and
$h_x(1)=s$.

\vskip 3mm \par\noindent {\bf Remark.} Strict convexity is needed
in (iii) and (iv). Indeed, let  $x=(0,2)$ and
$$K=conv\big(\{(1,1), (-1,1), (-2,0), (2,0)\}\big).$$ Then
$\KC=[(-1,1),(1,1)]$. However for any point $tx$ with $t>1$, $$
\overline{[tx, K]\setminus K}\cap
\partial K =\partial K\setminus [(-2,0),(2,0)] \supsetneq \KC. $$
Thus, for any function $f$ with $f>0$ on $[(-2,0),(-1,1)]$ and/or
$[(1,1),(2,0)]$, $h_x(\cdot)$ is not  right continuous on
$[1,\infty)$.

\vskip 2mm \par\noindent To see that strict convexity is needed
also in (iv), observe that $K^{f,1/12}=K$ in Example \ref{ex-1}.
Thus, for $x\in \partial K^{f,1/12}=\pt K$, we have $$
\m_f({\overline{[x, K]\setminus K}\cap
\partial K})=0 \neq \frac{1}{12}.$$

\section{Geometric interpretation of functionals on convex bodies}
We now give  geometric interpretations of functionals on convex
bodies, such as $L_p$ affine surface area and mixed $p$-affine
surface area
 for all $p \neq -n$ using the {\em non
convex} illumination surface bodies. While there are no geometric
interpretations for mixed $p$-affine surface area, many geometric
interpretations of $L_p$ affine surface area have been discovered
in the last years, all based on using {\em convex} bodies (e.g.,
\cite{MW1, SW4, SW5, WY}). The remarkable new fact here is that
now the bodies involved in the geometric interpretation are not
necessarily convex.  \bt \label{TH:0} Let $K$ be a  convex body in
$C^2_+$. Let $c>0$ be a constant, and $f:
\partial K \rightarrow \mathbb{R}$ be an integrable function such  that  $f \geq c$ $\m _K $-almost everywhere.  Then
 \begin{equation}\label{Form:0:0:0}
 \lim
_{s\rightarrow 0} c_n \frac{|K^{f,s}|-|K|}{
s^{\frac{2}{n-1}}}=\int_{\partial K} \frac{\k
_K(x)^{\frac{1}{n-1}}}{f(x)^{\frac{2}{n-1}}}\,d\m _K(x),
\end{equation} where $c_n=2 |B_2^{n-1}|^{\frac{2}{n-1}}$.  \et

\par \noindent {\bf Remark.} As $d\mu_K= f_K d \sigma$, we also have
 \begin{equation}\label{SURFACE:1:1}
 \lim
_{s\rightarrow 0} c_n \frac{|K^{f,s}|-|K|}{
s^{\frac{2}{n-1}}}=\int_{S^{n-1}}
\frac{f_K(u)^{\frac{n-2}{n-1}}}{f(N_K^{-1}(u))^{\frac{2}{n-1}}}\,d\s
(u) ,
\end{equation}
where
$N_K^{-1}$ is the inverse of the Gauss map $N_K(\cdot)$.

 \vskip 3mm
 \noindent
The geometric interpretation of $L_p$ affine surface area is then
a corollary to Theorem \ref{TH:0}. The theorem also gives
geometric interpretations of other known functionals on convex
bodies, e.g. the surface area and the mixed $p$-affine surface
area. Notice that  these geometric interpretations can also be
obtained using e.g. the  (convex) surface body \cite {SW5, WY}.
\vskip 3mm \noindent Define
$$
\tilde{f}(N^{-1}_K(u))=f_K(u)^{\frac{n-2}{2}}[f_p(K_1,u)\cdots
f_p(K_n,u)]^{\frac{1-n}{2(n+p)}},
$$
where $ f_p(K,u)= h_K(u)^{1-p}f_K(u)$.
\vskip 3mm
\bc\label{MIXAFF:1}
Let
$K$ and $ K_i$, $i=1, \cdots, n$,  be convex bodies in $C^2_+$.
Then
 \begin{equation*}\label{MIXAFF:1:1}
 \lim
_{s\rightarrow 0} c_n \frac{|K^{\tilde{f},s}|-|K|}{
s^{\frac{2}{n-1}}}=as_p(K_1, \cdots, K_n).
\end{equation*}
In particular, if all $K_i$ coincide with $K$, then $as_p(K_1, \cdots, K_n)= as_p(K)$ and
we get a geometric interpretation of $as_p(K)$
\begin{equation*}\label{cor1}
c_{n}\lim_{s \to 0} \frac{| K^{g_p, s}|-|K|} {s^\frac{2}{n-1}} =\
as_{p}(K),
\end{equation*}
where  $g_p: \partial K \rightarrow {\mathbb R}$ is defined by $g_p(x)=\k_K(x)^\frac{n+2p-np}{2(n+p)}\
  \langle x, N_K(x)\rangle^\frac{n(n-1)(p-1)}{2(n+p)}.$ \ec

\bc \label{SURFACE:0} Let $K$ be a convex body in $C^2_+$ and
$g(x)=\sqrt{\k_K(x)}$. Then
 \begin{equation*}\label{SURFACE:0:1}
 \lim
_{s\rightarrow 0} c_n \frac{|K^{g,s}|-|K|}{ s^{\frac{2}{n-1}}}=\m
_K(\pt K).
\end{equation*}
\ec

\noindent The proof of the corollaries follows immediately from
Theorem \ref{TH:0}. To prove Theorem \ref{TH:0}, we need several
other concepts and lemmas.

\vskip 2mm \noindent As  $K$ is in $ C^2_+$, for any $x \in
\partial K$, the indicatrix of Dupin is an ellipsoid. As in \cite{SW5}, we apply an affine transform $T:\mathbb
R^{n}\rightarrow\mathbb R^{n}$ to $K$ so that the indicatrix of
Dupin is transformed into an $(n-1)$-dimensional Euclidean ball.
$T$ has the following properties:
\begin{equation}\label{T-prop}
T(x)=x
\hskip 10mm
T(N_{K}(x))=N_{K}(x)
\hskip 10mm
\det(T)=1
\end{equation}
and  $T$ maps a measurable subset of a hyperplane orthogonal to
$N_{K}(x)$ onto a subset of the same $(n-1)$-dimensional measure.
It was also shown in \cite{SW5} that
for any $\epsilon>0$ there is $\Delta_1=\Delta_1(\varepsilon)>0$
such that for all measurable subsets $A$ of $\partial K\cap
H^{-}(x -\Delta_1 N_{ K} (x), N_{ K}(x))$
\begin{equation}\label{affineimage1}
(1-\epsilon)\  \mu_{K}(A) \leq |T(A)| \leq(1+\epsilon)\
\mu_{K}(A). \end{equation}
\par
$T(K)$ can be approximated at $x=T(x)$ by a $n$-dimensional
Euclidean ball: For any $\epsilon>0$ there is
$\Delta_2=\Delta_2(\varepsilon)$ such that
\begin{eqnarray}\label{capest}
&& B_{2}^{n}\left(x-r N_{ K}(x),r\right)
\cap H^{-}\left(x-\Delta_2 N_{ K}(x),
N_{ K}(x)\right) \nonumber   \\
&&  \hskip 10mm
\subseteq T(K)\cap H^{-}\left(x-\Delta_2 N_{ K}(x),
N_{K}(x)\right)    \\
&& \hskip 10mm
\subseteq B_{2}^{n}\left(x-R N_{ K}(x),R\right)
\cap H^{-}\left(x-\Delta_2 N_{ K}(x),
N_{ K}(x)\right), \nonumber
\end{eqnarray}
where $r=r(x)=\kappa_K(x)^{-\frac{1}{n-1}}$ and $ R=R(x)$ with $r
\leq R \leq (1+\epsilon) r$ . We put
\begin{equation}\label{D}
\Delta=\Delta(\varepsilon) = \mbox{min}\{\Delta_1,\Delta_2\}.
\end{equation}
Moreover,
for $x\in \pt K$, let
\begin{equation} \label{xs}
x_s\in \pt K^{f,s} \mbox{ be such that} \ \ x\in [0,x_s]\cap \pt K
\end{equation}
and define $\tilde{x}_s$ to be the orthogonal projection of $x_s$
onto the ray $\{y: y=x+t N_K(x), t\geq 0\}$. Clearly
$T(\tilde{x}_s)=\tilde{x}_s$, and the distance from $T(x_s)$ to
the hyperplane $H(x, N_K(x))$ is the same as the distance from
$x_s$ to this hyperplane.

\vskip 3mm

We say that a family of sets $E_s\subseteq \pt K$, $0<s\leq s_0$
shrinks nicely to a point $x\in \pt K$ (see \cite{FO}) if

\vskip 2mm  \par \noindent (SN1) $\mbox{diam} E_s \rightarrow 0$,
as $s \rightarrow 0$. \vskip 2mm \noindent (SN2) There is a
constant $\beta >0$ such that for all $s \leq s_0$ there exists
$t_s$ with
$$\mu_K\left(\pt K \cap
B(x,t_s)\right) \geq \mu_K(E_s) \geq \beta\  \mu_K\left(\pt K \cap
B(x,t_s)\right).$$

\par 
\noindent 
\bl \label{Le:4:1} Let $K$ be a convex
body in $C^2_+$ and $f: \partial K \rightarrow \mathbb{R}$ an
integrable, $\m _K $-almost everywhere strictly positive function.
Let  $x\in \pt K$ and let $x_s$ and $\tilde{x}_s$ be as above
(\ref{xs}). Then
\par
\noindent
(i)  The family $\pt K\cap
\overline{[\tilde{x}_s,K]\setminus K}$, $0 < s \leq s_0$  shrinks
nicely to $x$.
\par
\noindent
(ii)  The family
$\pt K\cap \overline{[x_s,K]\setminus K}$, $0 < s \leq s_0$
shrinks nicely to $x$.
\par
\noindent
(iii)
\begin{equation}\label{Form:0:2}
\lim _{s\rightarrow 0}\frac{\m_f(\pt K\cap
\overline{[\tilde{x}_s,K]\setminus K})}{\mu_K(\pt K\cap
\overline{[\tilde{x}_s,K]\setminus K})}=f(x) \quad \m_K \mbox{-a.e.}
\end{equation}
\par
\noindent
(iv)
\begin{equation}\label{Form:0:2a}
\lim _{s\rightarrow 0}\frac{\m_f(\pt K\cap
\overline{[x_s,K]\setminus K})}{\mu_K(\pt K\cap
\overline{[x_s,K]\setminus K})}=f(x) \quad \m_K \mbox{- a.e.}
\end{equation} \el

\vskip 2mm \noindent {\bf Proof.}
Formulas (\ref{Form:0:2}) and (\ref{Form:0:2a}) in
(iii) and (iv) follow  from the Lebesgue differentiation theorem
(see \cite{FO}) once  we  have proved that $\pt K\cap
\overline{[\tilde{x_s},K]\setminus K}$ and $\pt K\cap
\overline{[x_s,K]\setminus K}$  shrink nicely to $x$. Therefore it
is enough to prove (i) and (ii).

\vskip 2mm \par \noindent
(i) For $x\in \pt K$, let $r=r(x)$ and
$R=R(x)$ be as in (\ref{capest}). We abbreviate $B(r)=
B_{2}^{n}(x-r N_{ K}(x), r)$ and $B(R) = B_{2}^{n}(x-R N_{ K}(x),
R)$. Let
\begin{equation*}
\D(x,s)=\left\langle \frac{x}{\|x\|}, N_{ K}(x)\right\rangle
\|x_s-x\|=\langle x_s-x,N_K(x)\rangle
\end{equation*}
be the distance from $x_s$ to $H\big(x, N_{ K}(x)\big)$. This is
the same as the distance from $\tilde{x}_s$ (defined after formula
(\ref{xs})) to $H\big(x, N_{ K}(x)\big)$.
\par\noindent Let $ h_R = h_R(s) = \frac{R\ \D(x,s)}{R+\D(x,s)} $ be the height
of the cap of $B(R)$ that is ``illuminated" by $\tilde{x}_s$. Then
\begin{equation}\label{F:1}
H^{-}\big(x-h_R N_{ K}(x), N_{ K}(x)\big)\cap \pt B(R)=\overline{[T(\tilde{x}_s),
B(R)]\setminus B(R)} \cap \pt B(R).
\end{equation}
\par \noindent
Let $\Delta$ be as in (\ref{D}). Since $\D(x,s)\rightarrow 0$ as
$s\rightarrow 0 $, one can choose $s _1\leq s_0$, such that for
all $0 < s \leq s _1$, $h=2 h_R < \Delta$. Therefore
(\ref{capest}) holds:
\begin{eqnarray}
& & H^{-}(x-h N_{ K}(x), N_{ K}(x))\cap B(r)\subset  H^{-}\big(x-h
N_{ K}(x), N_{ K}(x)\big)\cap T(K)\nonumber \\ & &\qquad \ \ \ \ \
\ \  \subset H^{-}\big(x-h N_{ K}(x), N_{ K}(x)\big)\cap B(R).
\label{F:0}
\end{eqnarray}

\vskip 2mm \noindent (\ref{F:0}) and (\ref{F:1}) imply that
for all small enough $s \leq s_2 \leq s_1$
$$
T\big(\overline{[\tilde{x}_s, K] \setminus K}\cap \pt
K\big)=\overline{[T(\tilde{x}_s), T(K)] \setminus T(K)}\cap \pt
T(K) $$
$$\subseteq H^{-}\big(x-h N_{ K}(x), N_{ K}(x)\big) \cap \partial T(K)
\subseteq B(R) \cap H^{-}\big(x-h N_{ K}(x), N_{ K}(x)\big).
$$
\vskip 2mm \par \noindent Let $t_s=\|x-z\|=\sqrt{\frac{4R^2 \
\D(x,s)}{R+\D (x,s)}}$ where $z$ is any point in $H\big(x-h N_{
K}(x), N_{ K}(x)\big)\cap \partial B(R)$. As $B(R) \cap
H^{-}\big(x-h N_{ K}(x), N_{ K}(x)\big) \subseteq B^n_2(x,t_s)$,
$$ T\big(\overline{[\tilde{x}_s, K]\setminus K} \big) \cap \pt
T(K) = T\big(\overline{[\tilde{x}_s, K]\setminus K}\cap \pt K\big
)\subseteq B^n_2(x,t_s) \cap \partial T(K)
$$ and $t_s \rightarrow 0$ as $s \rightarrow 0$.
\newline
This shows that condition (SN1)
is satisfied for  $T\big(\overline{[\tilde{x}_s, K]\setminus K}
\big) \cap \partial T(K)$  to shrink nicely to $T(x)=x$.

\vskip 2mm \par \noindent
We now show that condition (SN2) also
holds true.
\par \noindent
First, $\overline{[T(\tilde{x}_s),
B(r)]\setminus B(r)}\cap H(x,N_K(x))$ is a $(n-1)$-dimensional Euclidean ball
with radius $\frac{r \D (x,s)}{\sqrt{2r\D(x,s)+\D^2(x,s)}}$. Then
for any $0<s<s_2$,
\begin{eqnarray*} |T(\overline{[\tilde{x}_s, K]\setminus K}\cap \pt K)|
&\geq& |\overline{[T(\tilde{x}_s), B(r)]\setminus B(r)}\cap
H(x,N_K(x))|\\ &\geq& |B_2^{n-1}| \left(\frac{r \D
(x,s)}{\sqrt{2r\D (x,s)+\D ^2 (x,s)}}\right)^{n-1}.
\end{eqnarray*}
\par\noindent
We can choose ( a new, smaller) $s_2$  such that $\frac{\D(x,s)}{r}\leq
2$.  Then
\begin{eqnarray} \label{lemma:below}
|T(\overline{[\tilde{x}_s,
K]\setminus K}\cap \pt K)|\geq 2^{1-n} \ |B_2^{n-1}| \
(r\D(x,s))^{\frac{n-1}{2}}.
\end{eqnarray}
\par \noindent
On the other hand, for $\varepsilon$ small enough, there exists $s_3<s_2$, such
that, for all  $0<s\leq s_3$ and  for any subset $A$ of
$H^{-}(x-2hN_K(x), N_K(x))\cap \pt T(K)$ \cite{SW5}
\begin{equation}\label{Proj-1}
|P_{H(x-2hN_K(x), N_K(x))}(A)|\leq |A| \leq
(1+\varepsilon)|P_{H(x-2hN_K(x), N_K(x))}(A)|\end{equation} where
$P_H(A)$ is the orthogonal projection of $A$ onto the hyperplane $H$.
We apply this to $A=B^n_2(x,t_s)\cap \partial T(K)$:
\begin{eqnarray}
|B^n_2(x,t_s)\cap \partial T(K)|&\leq
&(1+\varepsilon)|P_{H(x-2hN_K(x), N_K(x))}(B^n_2(x,t_s)\cap
\partial T(K))|\nonumber \\&\leq&
(1+\varepsilon) \ |B(R)\cap H(x-2hN_K(x),N)|\nonumber \\
&\leq&(1+\varepsilon )\ |B_2^{n-1}| \ \left(
\frac{2\sqrt{2}R \sqrt{R\D (x,s)}}{R+\D(x,s)}\right)^{n-1}\nonumber \\
&\leq & 4^n \ |B_2^{n-1}|\ (r\D(x,s))^{\frac{n-1}{2}}
\label{F:1:1}\end{eqnarray}
The last inequality follows as
$r \leq R<(1+\varepsilon)r$. It now follows from (\ref{lemma:below}) and (\ref{F:1:1}) that
also  condition (SN2) holds true for e.g. $\beta=8^{-n}$.

\par \noindent
Hence the family $T\big(\overline{[\tilde{x}_s, K]\setminus K}
\big) \cap \partial T(K)$ shrinks nicely to $T(x)=x$ and
therefore, as $T^{-1}$ exists, the family $\overline{[\tilde{x}_s,
K]\setminus K} \cap \partial K = T^{-1} \big(
T\big(\overline{[\tilde{x}_s, K]\setminus K} \big) \cap \partial
T(K) \big) $ shrinks nicely to $x$.

\vskip 3mm \noindent (ii) Let  $v_1=x_s-(x-r N_{K}(x))$ and
$v_2=x_s-( x-RN_{K}(x))$. $\theta$  denotes the angle between
$N_K(x)$ and $x$ and $\phi_i=\phi_i(x,s), i=1,2$ is the angle
between  $N_K(x)$ and $v_i,  i=1,2$. These angles  can be computed
as follows
\begin{eqnarray*}
\mbox{tan}(\phi_1)&=&\frac{\D (x,s) \ \mbox{tan}(\theta)}{r+\D(x,s)}  \\
\mbox{tan}(\phi_2)&=&\frac{\D (x,s) \ \mbox{tan}(\theta)}{R+\D(x,s)}.
\end{eqnarray*}
Then, for $i=1,2$,  $\phi _i \rightarrow 0$ as $s\rightarrow 0$.
Since $K$ is in $C_+^2$, this means that  for any $\varepsilon >0$
there is $\bar{s}_\varepsilon \leq s_0$ such that for all $s\leq
\bar{s}_\varepsilon$
\begin{equation}\label{XS-XST-1}
1-\varepsilon \leq \frac{\mu_K(\overline{[x_s,K]\setminus K}\cap
\partial K)}{\mu_K(\overline{[\tilde{x}_s,K]\setminus K}\cap
\partial K)} \leq 1+\varepsilon.
\end{equation}
By (\ref{capest}) and as $\phi _i \rightarrow 0, i=1,2$ as $s\rightarrow 0$,
one can choose $\tilde{h}_R = \tilde{h}_R(s) = \frac{3R\
\D(x,s)}{R+\D(x,s)}$ so small that
\begin{eqnarray*}
T\left(\overline{[{x}_s, K] \setminus K}\cap \pt
K\right)&=&\overline{[T({x}_s), T(K)]\setminus T(K)}\cap \pt
T(K)\\
&\subset& H^{-}\left(x-\tilde{h}_R  N_K(x), N_K(x)\right)\cap
B(R). \end{eqnarray*} Let $ \tilde{t}_s = \sqrt{\frac{6 R^2 \
\D(x,s)}{R+\D (x,s)}} $ be the distance from $x$ to any point in
$H\left(x - \tilde{h}_R N_K(x), N_K(x) \right) \cap \partial
B(R)$. Then
$$T\left(\overline{[{x}_s, K]\setminus K} \cap \pt K\right)\subseteq
B^n_2(x, \tilde{t}_s) \cap \partial T(K). $$ \noindent
(\ref{affineimage1}),  (\ref{XS-XST-1}) and Lemma \ref{Le:4:1} (i)
then give
\begin{eqnarray}|T\left(\overline{[{x}_s, K]\setminus K} \cap \pt
K\right)|&\geq& (1-\varepsilon)^3 |T\left(\overline{[\tilde{x}_s,
K]\setminus K} \cap \pt K\right)|\nonumber \\ &\geq &
(1-\varepsilon)^3 \beta |B^n_2(x,t_s)\cap
\partial T(K)| \label{F:1:1:3}.
\end{eqnarray}
Furthermore, by (\ref{Proj-1}), one has \begin{eqnarray*}
|B^n_2(x,t_s)\cap \pt T(K)|&\geq& |H(x-hN_K(x), N_K(x))\cap
T(K)|\\ &\geq & |H(x-hN_K(x), N_K(x))\cap B(r)|\\
&=&\left(\frac{4R^2
r\D(x,s)+4Rr\D(x,s)^2-4R^2\D(x,s)^2}{(R+\D(x,s))^2}\right)^{\frac{n-1}{2}}|B^{n-1}_2|.
\end{eqnarray*}
Since $\D(x,s)\rightarrow 0$ as $s\rightarrow 0$, for
$\bar{s}_{\varepsilon}$ small enough, and any
$0<s<\bar{s}_{\varepsilon}$, one get
\begin{eqnarray}
|B^n_2(x,t_s)\cap \pt T(K)|&\geq&
\left(\frac{2R}{R+\D(x,s)}\sqrt{r\D(x,s)-\D(x,s)^2}\right)^{{n-1}}|B^{n-1}_2|\nonumber\\
&\geq &
2^{-n}\left(r\D(x,s)\right)^{\frac{n-1}{2}}|B^{n-1}_2|\label{F:1:1:2}.
\end{eqnarray}
A computation similar to (\ref{F:1:1}) shows that for all $0<s\leq
\bar{s}_{\varepsilon}$ with (a possibly new)
$\bar{s}_{\varepsilon}$ small enough
\begin{eqnarray} |B^n_2(x,\tilde{t}_s)\cap \partial T(K)|&\leq&
(1+\varepsilon) \ |B(R)\cap H(x-\tilde{h}N_K(x),N)|\nonumber \\
&=&(1+\varepsilon )\ |B_2^{n-1}| \ \left(
\frac{\sqrt{6} R \sqrt{R\D (x,s)}}{R+\D(x,s)}\right)^{n-1}\nonumber \\
&\leq & 3^n \ |B_2^{n-1}|\ (r\D(x,s))^{\frac{n-1}{2}}
\label{F:1:1:1}\end{eqnarray}
\par\noindent
 (\ref{F:1:1:3}), (\ref{F:1:1:2}) and (\ref{F:1:1:1})
imply  that
\begin{equation*}|T\left(\overline{[{x}_s, K]\setminus
K}\cap \pt K\right)| \geq (48)^{-n-1}\beta |B^n_2(x,\tilde{t}_s)\cap
\partial T(K)|.
\end{equation*}
This shows that  $T\left(\overline{[{x}_s, K]\setminus
K}\cap \pt K\right)$ shrinks  nicely to $x$. Therefore  also $\overline{[{x}_s, K]\setminus
K}\cap \pt K$ shrinks nicely to $x$.

\vskip 3mm
 \bl\label{Le:5}
Let $K$ be a convex body  in $C^2_+$ and $f: \partial K
\rightarrow \mathbb{R}$ an integrable, $\m _K $-almost everywhere
strictly positive function. Then for  $\mu_K$-almost all $x\in \pt
K$ one has
\begin{equation}\label{Form:0:3}\lim
_{s\rightarrow 0} c_n \frac{\langle x,
N_{K}(x)\rangle\left[\left(\frac{\|x_s\|}{\|x\|}\right)^n-1\right]}{
s^{\frac{2}{n-1}}}=\frac{\k
_K(x)^{\frac{1}{n-1}}}{f(x)^{\frac{2}{n-1}}},
\end{equation} where $x_s \in \partial K^{f,s}$ is such that $x \in [0,x_s]$. \el
\vskip 2mm
\noindent
{\bf Proof.} It is enough to consider  $x\in \pt K$ such that
$f(x)>0$.
As $x$ and $x_{s}$ are collinear and as
$(1+t)^n\geq 1+tn$ for $t\in [0,1)$, one has for small enough $s$,
\begin{eqnarray*}
\frac{\langle x,
N_{K}(x)\rangle}{n}\left[\left(\frac{\|x_s\|}{\|x\|}\right)^n-1\right]=\frac{\langle
x,
N_{K}(x)\rangle}{n}\left[\left(1+\frac{\|x_s-x\|}{\|x\|}\right)^n-1\right]\geq
\D(x,s) .
\end{eqnarray*}
Recall that $\D(x,s)=\left\langle \frac{x}{\|x\|}, N_{
K}(x)\right\rangle \|x_s-x\|=\langle x_s-x,N_K(x)\rangle$ is the
distance from $x_s$ to $H\left(x,N_{ K}(x)\right)$.
\par\noindent
Similarly, as $(1+t)^n \leq 1+nt+2^n t^2$ for $t\in [0,1)$, one has for $s$ small enough,
\begin{equation}\label{Form:0:0}
\frac{\langle x,
N_{K}(x)\rangle}{n}\left[\left(\frac{\|x_s\|}{\|x\|}\right)^n-1\right]\leq
\D(x,s) \left[ 1+\frac{2^n}{n} \ \left(\frac{\|x_s
-x\|}{\|x\|}\right) \right].\
\end{equation}
\par\noindent
Hence for  $\varepsilon>0$ there exists $s_\varepsilon \leq s_0$ such that for all
$0<s\leq s_\varepsilon$
\begin{equation*}
1\leq \frac{\langle x,
N_{K}(x)\rangle\left[\left(\frac{\|x_s\|}{\|x\|}\right)^n-1\right]}{n\D(x,s)}\leq
1+\varepsilon.
\end{equation*}
$K$ is strictly convex as $K \in C^2_+$. Thus,
$\mu_f(\partial K \cap \overline{[x_s, K] \setminus K} ) =s$ by Lemma \ref{Lemma:prop} (iv). Therefore
$$
1\leq \frac{\langle x,
N_{K}(x)\rangle\left[\left(\frac{\|x_s\|}{\|x\|}\right)^n-1\right]
\left(\mu_f(\partial K \cap \overline{[x_s, K] \setminus
K})\right)^\frac{2}{n-1}}{n\ s^\frac{2}{n-1}\ \D(x,s)}\leq
1+\varepsilon.
$$
By Lemma \ref{Le:4:1} (iv) and  (\ref{XS-XST-1}), it then
follows that we can choose (a new) $s_\varepsilon$ so small that we have
for all $s \leq s_\varepsilon$
\begin{equation}\label{Form-1} 1-c_1 \varepsilon \leq \frac{\langle x,
N_{K}(x)\rangle\left[\left(\frac{\|x_s\|}{\|x\|}\right)^n-1\right]
\left(f(x)\ \mu_K(\partial K \cap \overline{[\tilde{x_s}, K]
\setminus K})\right)^\frac{2}{n-1}}{n\ s^\frac{2}{n-1}\
\D(x,s)}\leq 1+c_2 \varepsilon. \end{equation} with absolute
constants  $c_1$ and $c_2$.

\par \noindent Let $T$ be as in (\ref{T-prop}) and let
$r=r(x)$ and $R=R(x)$ be as in (\ref{capest}). We abbreviate again
$B(r)= B_{2}^{n}(x-r N_{ K}(x), r)$ and $B(R) = B_{2}^{n}(x-R N_{
K}(x), R)$. Let $h_r = h_r(s) = \frac{r\D(x,s)}{r+\D(x,s)}$. As
$h_r \rightarrow 0$ as $s\rightarrow 0$, we have  for all $s$
sufficiently small that  $h_r < \Delta $ where $\Delta $ is as in
(\ref{D}). Hence by (\ref{capest})
$$
H^-\left(x - h_r N_K(x), N_K(x)\right )\cap  B(r) \subset
H^-\left(x - h_r N_K(x), N_K(x)\right ) \cap T(K).
$$
If we denote by $P_H$ the orthogonal projection onto the hyperplane
$H^-\left(x - h_r N_K(x), N_K(x)\right )$, this then implies that
\begin{eqnarray*}
P_{H} \big(\partial K \cap \overline{[\tilde{x}_s,K] \setminus K}
\big) &\supset & P_{H} \big(\partial T^{-1} (B(r)) \cap
\overline{[\tilde{x}_s,T^{-1}(B(r))] \setminus T^{-1}(B(r))}\big)
\\&=& H\left(x - h_r N_K(x), N_K(x)\right )\  \cap \ T^{-1}(B(r)).
\end{eqnarray*}
Hence for $s$ sufficiently small
\begin{eqnarray}
\m_K(\partial K\cap \overline{[\tilde{x}_s,K]\setminus
K}) &\geq& \big|P_H \left(\partial K\cap \overline{[\tilde{x}_s,K]\setminus
K}\right) \big| \nonumber \\
&\geq& \big |H\left(x - h_r N_K(x), N_K(x)\right )\  \cap \ T^{-1}(B(r)) \big|\nonumber \\
&=&\big |T\big(H\left(x - h_r N_K(x), N_K(x)\right )\  \cap \ T^{-1}(B(r))\big) \big| \nonumber \\
&=& \big | H\left(x - h_r N_K(x), N_K(x)\right )\  \cap B(r)\big| \nonumber \\
&=&\frac{\big(1-\frac{\D(x,s)}{2 r}\big)^\frac{n-1}{2}}{\big(1+\frac{\D(x,s)}{ r}\big)^{n-1}}\
\left (2 r \D(x,s)\right)^{\frac{n-1}{2}}        |B^{n-1}_2| \nonumber \\
&\geq &(1-c_3 \varepsilon) \ \frac {\left(2
\D(x,s)\right)^{\frac{n-1}{2}} } {\sqrt{\k _K(x)} }    \
|B^{n-1}_2|     \label{K-3}
\end{eqnarray}
where $c_3>0$ is an absolute constant.

\vskip 3mm \par\noindent By (\ref{capest}), one has
\begin{eqnarray*}
P_{H}\big(\partial K\cap \overline{[\tilde{x}_s,K]\setminus
K}\big) &\subseteq & P_{H} \big(\partial \left(T^{-1}(B(R))\right)
\cap \overline{[\tilde{x}_s, T^{-1}(B(R))]\setminus
T^{-1}(B(R))}\big)\\ &=& H\left (x-h_R N_K(x), N_K(x)\right) \cap
T^{-1}(B(R)) \end{eqnarray*} where  $H=H\left (x-h_R N_K(x),
N_K(x)\right)$. The equality follows as
$h_R=\frac{R\D(x,s)}{R+\D(x,s)}$.

\vskip 2mm \par\noindent Together with (\ref{Proj-1}), for
$s_\varepsilon$ small enough, whenever $0<s<s_\varepsilon$, one
has
\begin{eqnarray*}
\mu_K(\partial K\cap \overline{[\tilde{x}_s,K] \setminus K})&\leq
& (1+\varepsilon) |P_{H}(\partial K\cap \overline{[\tilde{x}_s,K]
\setminus K}) | \nonumber \\ &\leq & (1+\varepsilon) \left |
H\left (x-h_R N_K(x), N_K(x)\right) \cap T^{-1}(B(R))\right|.
\end{eqnarray*}

\par\noindent
A calculation similar to (\ref{K-3}) then shows that with an
absolute constant $c_4$
\begin{equation}\label{K-5}
\m_K(\partial K\cap \overline{[\tilde{x}_s,K]\setminus K}) \leq
(1+c_4 \varepsilon)  \ \frac{ (2\D(x,s))^{\frac{n-1}{2}}}
{\sqrt{\k _K(x)}} \  |B_2^{n-1}|.
\end{equation}
\par\noindent
Combining (\ref{Form-1}), (\ref{K-3}) and (\ref{K-5}), we prove
the formula (\ref{Form:0:3}), i.e.,
$$\lim _{s\rightarrow 0} c_n \frac{\langle x,
N_{K}(x)\rangle\left[\left(\frac{\|x_s\|}{\|x\|}\right)^n-1\right]}{n\
s^{\frac{2}{n-1}}}=\frac{\k
_K(x)^{\frac{1}{n-1}}}{f(x)^{\frac{2}{n-1}}}.$$

\vskip 5mm
 \bl\label{Le:6}
Let $K$ be a convex body in $C^2_+$. Let $c>0$ be a constant and
$f: \partial K \rightarrow \mathbb{R}$ an integrable function with
$f\geq c$ $\m _K $-almost everywhere. Then there exists
$\overline{s} \leq s_0$, such that for all $s\leq \overline{s}$,
\begin{equation*}
 \frac{\langle x,
N_{K}(x)\rangle\left[\left(\frac{\|x_s\|}{\|x\|}\right)^n-1\right]}{
s^{\frac{2}{n-1}}} \leq c(K,n),
\end{equation*}
where $c(K,n)$ is a constant (depending on $K$ and $n$ only), and
$x$ and $x_{s}$ are as in Lemma \ref{Le:5}. \el
\vskip 2mm
\noindent
{\bf Proof.}

As $K\in C^2_+$, by the Blaschke rolling theorem \cite{Sch}, there
exists $r_0>0$ such that for all $x \in \partial K$,
$ B^n_2(x - r_0 N_K(x), r_0) \subseteq K. $ Let $\gamma$ be such
that  $0<\g \leq \mbox{min}\{1,r_0\}$. By Lemmas \ref{Le:1} (i)
and \ref{Le:2} (ii), $K=K^{f,0}=\bigcap _{s>0}K^{f,s}$. Therefore
there exists $\overline{s} = s_\gamma \leq s_0$, such that for all
$s \leq \overline{s}$, $K^{f,s}\subseteq (1+\g)K$. Hence for
$x_s\in \partial K^{f,s}$ and $x=[0,x_s]\cap \partial K$, $
\frac{\|x_s\|}{\|x\|}\leq 1+\g $, or equivalently -as $x$ and
$x_s$ are collinear-
\begin{equation}\label{gamma2}
\frac{\|x_s\|}{\|x\|} - 1 = \frac{\|x_s-x\|}{\|x\|} \leq \gamma \leq 1.
\end{equation}
Together with (\ref{Form:0:0}), one has for all $s \leq
\overline{s}$ (with a possibly smaller $\overline{s}$)
\begin{equation}\label{gamma}
\langle x,
N_{K}(x)\rangle \left[\left(\frac{\|x_s\|}{\|x\|}\right)^n-1\right]\leq
\D(x,s) \left[ n+2^n \right].\
\end{equation}
\par
As $K \in C^2_+$, $K$ is strictly convex. Hence, by Lemma
\ref{Lemma:prop} (iv),  (\ref{capest}) and  as $f\geq c$ on $\partial
K$
\begin{eqnarray*}
s &=& \m_f(\partial K \cap \overline{[x_s,K]\setminus  K}) = \int_{ \partial K \overline{[x_s,K]\setminus K})} \ f d\mu_K
\geq  c\ \mu_K \big(\partial K \cap \overline{[x_s,K]\setminus K}\big)\\
&\geq & c\  \big|H\left(x, N_K(x) \right)\cap [x_s, \
T^{-1}(B(r))] \big| = c\  \big|H\left(x, N_K(x) \right)\cap
[T(x_s), \ B(r)] \big|
\end{eqnarray*}
where $T$ is as in  (\ref{T-prop}) and $r=r(x)$ is as in
(\ref{capest}). \vskip 2mm \par\noindent As in the proof of Lemma
\ref{Le:2} (i), $H\left(x, N_K(x) \right)\cap [T(x_s), \ B(r)] $
contains a $(n-1)$-dimensional Euclidean ball of radius at least
$$
\frac{r\sqrt{2r\D(x,s)+\D^2(x,s)}}{2r+\D(x,s)}\geq
\frac{\a}{1+2\a}\sqrt{2r_0\D(x,s)}.
$$
The inequality follows as $
\D(x,s)=\frac{\|x_s-x\|}{\|x\|} \langle x, N_K(x) \rangle\leq
\frac{\g}{\a}\leq \frac{1}{\a},$ which is a direct consequence of
(\ref{alpha}) and (\ref{gamma2}).

\vskip 2mm  \par \noindent Hence $ s \geq
c\left(\frac{\a}{1+2\a}\sqrt{2r_0\D(x,s)}\right)^{n-1} |B^{n-1}_2|
$ and with (\ref{gamma}) we get that
$$
\frac{\langle x,
N_{K}(x)\rangle}{s^\frac{2}{n-1}}\  \left[\left(\frac{\|x_s\|}{\|x\|}\right)^n-1\right]\leq
(n+2^n) \left(\frac{1+ 2 \alpha}{\alpha}\right)^2 \big( 2 r_0\  c^\frac{2}{n-1} |B^{n-1}_2| ^\frac{2}{n-1} \big)^{-1}.
$$

\vskip 4mm
\noindent
Finally, we also need the following lemma. It is well known and we omit the proof.
\vskip 3mm
\noindent
\bl \label{Le:4}
Let $K$ be a convex body and $L$ be a star-convex body in
 $\bbR^n$.\\
 (i) If $0\in int(L)$ and $L\subset K$, then
 $$
 |K|-|L|=\frac{1}{n}\int _{\partial K} \langle x,
N_{K}(x)\rangle\left[1-\left(\frac{\|x'\|}{\|x\|}\right)^n\right]\,d\m_{K}
(x)
$$
where $x\in \partial K$ and  $x'\in \partial L\cap [0,x]$. \\
(ii) If $0\in int(K)$ and $K\subset L$, then
$$
|L|-|K|=\frac{1}{n}\int _{\partial K} \langle x,
N_{K}(x)\rangle\left[\left(\frac{\|x'\|}{\|x\|}\right)^n-1\right]\,d\m_{K}
(x)
$$
where $x\in \partial K$, $x'\in \partial L$ and  $ x=\partial
K\cap [0,x']$. \el

\par \noindent {\bf Proof of Theorem \ref{TH:0}.}
\par \noindent As $K \in C^2_+$, $K$ is strictly convex. By Lemmas
\ref{Le:4}, \ref{Le:5}, \ref{Le:6} and Lebegue's  Dominated
Convergence theorem
\begin{eqnarray*}
c_n \lim _{s\rightarrow 0
}\frac{|K^{f,s}|-|K|}{s^{\frac{2}{n-1}}}&=&c_n \lim _{s\rightarrow
0 } \ \int _{\partial K} \frac{\langle x,
N_{K}(x)\rangle\left[\left(\frac{\|x_s\|}{\|x\|}\right)^n-1\right]}{n\
s^{\frac{2}{n-1}}}\,d\m_{K} (x)\\&=&c_n \ \int _{\partial K} \lim
_{s\rightarrow 0 } \frac{\langle x,
N_{K}(x)\rangle\left[\left(\frac{\|x_s\|}{\|x\|}\right)^n-1\right]}{n\
s^{\frac{2}{n-1}}}\,d\m_{K} (x)\\ &=&\int _{\partial K}\frac{\k
_K(x)^{\frac{1}{n-1}}}{f(x)^{\frac{2}{n-1}}}\, d \m_K(x).
\end{eqnarray*}

\vskip 2mm \noindent Elisabeth Werner\\
{\small Department of Mathematics \ \ \ \ \ \ \ \ \ \ \ \ \ \ \ \ \ \ \ Universit\'{e} de Lille 1}\\
{\small Case Western Reserve University \ \ \ \ \ \ \ \ \ \ \ \ \ UFR de Math\'{e}matique }\\
{\small Cleveland, Ohio 44106, U. S. A. \ \ \ \ \ \ \ \ \ \ \ \ \ \ \ 59655 Villeneuve d'Ascq, France}\\
{\small \tt elisabeth.werner@case.edu}\\ \\
\and Deping Ye\\
{\small Department of Mathematics}\\
{\small Case Western Reserve University}\\
{\small Cleveland, Ohio 44106, U. S. A.}\\
{\small \tt dxy23@case.edu}
\end{document}